\documentclass{amsart}
\usepackage{amsmath}
\usepackage{amsthm}
\usepackage[applemac]{inputenc} 
\usepackage{hyperref}
\usepackage{nicematrix}
\usepackage{amssymb}
\usepackage{graphicx}
\usepackage{enumitem}
\usepackage{amsfonts}
\usepackage[all]{xy}
\usepackage{array}
\usepackage{tikz-cd}
\usepackage{yhmath}
\usepackage{graphicx}
\usepackage{xcolor}
\usepackage{ytableau}

\usepackage[textsize=tiny]{todonotes}

\usepackage{enumitem}

\newcommand\lra{\longrightarrow}

\newcommand{\ol}{\overline}

\newcommand{\C}{\mathbb{C}}

\renewcommand{\P}{\mathbb{P}}

\newcommand{\Q}{\mathbb{Q}}
\newcommand{\R}{\mathbb{R}}

\newcommand{\Z}{\mathbb{Z}}

\newcommand{\cA}{\mathcal A}

\newcommand{\cG}{\mathcal G}

\newcommand{\cO}{\mathcal O}

\newcommand{\cY}{\mathcal Y}

\newcommand{\fg}{\mathfrak g}

\newcommand{\fn}{\mathfrak n}

\newcommand{\DA}{{\rm A}}
\newcommand{\DB}{{\rm B}}
\newcommand{\DC}{{\rm C}}

\newcommand{\DG}{{\rm G}}

\newcommand{\Ad}{\operatorname{Ad}}

\DeclareMathOperator{\Eff}{Eff}

\DeclareMathOperator{\Mov}{Mov}
\DeclareMathOperator{\mov}{mov}

\DeclareMathOperator{\Nef}{Nef}

\DeclareMathOperator{\NU}{N^1}
\DeclareMathOperator{\Pic}{Pic}

\newcommand{\CX}{\mathcal{C}\mkern-1mu X}

\newcommand{\GU}{\mathcal{G}\mkern-1mu U}
\newcommand{\GX}{\mathcal{G}\mkern-1mu X}
\newcommand{\MX}{\mathcal{M}\mkern-1mu X}

\newcommand{\MU}{\mathcal{M}\mkern-1mu U}

\newcommand{\LX}{\mathcal{L}\mkern-1mu X}

\renewcommand{\ss}{\mathrm{ss}}

\DeclareMathOperator{\Cl}{Cl}

\DeclareMathOperator{\codim}{codim}

\DeclareMathOperator{\Cox}{Cox}

\DeclareMathOperator{\HH}{H}

\DeclareMathOperator{\Hom}{Hom}
\DeclareMathOperator{\id}{id}

\DeclareMathOperator{\Int}{Int}

\DeclareMathOperator{\Proj}{Proj}

\DeclareMathOperator{\rk}{rk}

\DeclareMathOperator{\Spec}{Spec}

\DeclareMathOperator{\Mo}{M}
\DeclareMathOperator{\No}{N}

\DeclareMathOperator{\Env}{Env}

\newcolumntype{C}{>{$}c<{$}}

\newcommand\marit{11}
\newcommand\maren{18}

\newtheorem{theorem}{Theorem}[section]

\newtheorem{setup}[theorem]{Setup}
\newtheorem{lemma}[theorem]{Lemma}
\newtheorem{proposition}[theorem]{Proposition}
\newtheorem{corollary}[theorem]{Corollary}

\theoremstyle{definition}
\newtheorem{definition}[theorem]{Definition}
\newtheorem{example}[theorem]{Example}
\newtheorem{remark}[theorem]{Remark}
\newtheorem{notation}[theorem]{Notation}

\title{Maximal geometric quotients of flag varieties}

\author[Barban]{Lorenzo Barban}
\address{Center for Complex Geometry, Institute for Basic Science (IBS), 55 Expo-ro, Yuseong-gu, Daejeon, 34126, Republic of Korea}
\email{lorenzobarban@ibs.re.kr}

\author[Bianco]{Michele Bianco}
\address{Dipartimento di Matematica - Universit\'a di Trento, via Sommarive 14 I-38123  Trento (TN), Italy}
\email{michele.bianco@unitn.it}

\author[Sol\'a Conde]{Luis E. Sol\'a Conde}
\address{Dipartimento di Matematica - Universit\'a di Trento, via Sommarive 14 I-38123  Trento (TN), Italy}
\email{eduardo.solaconde@unitn.it}

\subjclass[2020]{Primary 14L30, 14E30; Secondary 14J30, 14L24, 14M17}
\thanks{The first author was supported by the Institute for Basic Science (IBS-R032-D1). Second and third author partially supported by INdAM--GNSAGA}

\begin{document}

\begin{abstract} 
	We introduce the notion of maximal geometric quotient, namely a geometric quotient of a polarized pair under a torus action which has a birational contraction to any other GIT quotient. For the action of the maximal torus on a complete flag variety, we show that a maximal geometric quotient exists and it is a Mori dream space, and its Cox ring is isomorphic to the Cox ring of the flag variety. We discuss birational properties of the maximal geometric quotients, their relations with the Chow quotients and explain how to construct them for maximal torus actions on flag varieties.
\end{abstract}
\maketitle

\tableofcontents

\section{Introduction}\label{sec:intro}

Geometric Invariant Theory (GIT), developed by Mumford in the 1960s (cf. \cite{MFK}), has played, over the years, a major role both in the study of the  birational geometry of algebraic varieties (\cite{Thaddeus1996}, \cite{Wlodarczyk}) and in the construction of moduli spaces (\cite{MFK}). However, GIT quotients depend on the choice of a linearization, and different linearizations may  produce different, yet birationally equivalent, GIT quotients. With this in mind, Kapranov, Sturmfels and Zelevinsky (cf. \cite{KSZ}) introduced a notion of quotient that does not depend on this choice, namely the {\em Chow quotient}. Intuitively, given an algebraic group action on a projective variety, the Chow quotient parametrizes the closures of general orbits and their limits in the Chow scheme of the variety.  Moreover, the Chow quotient dominates all GIT quotients (cf. \cite[Theorem~0.4.3]{Kapranov}), and for torus actions on projective varieties the Chow quotient is isomorphic to the limit quotient (cf. \cite[Corollary~2.6]{BakerChow}). Over the years, several Chow quotients have been computed: for instance, Kapranov, Sturmfels and Zelevinsky (cf. \cite{KSZ}) have studied the Chow quotients of projective toric varieties by torus actions; Kapranov (cf. \cite{Kapranov}) has related the Chow quotient of the Grassmannian $G(2,n)$, which parametrizes $2$-dimensional linear subspaces of $\C^n$, under the action of $(\C^*)^{n-1}$ with the Grothendieck--Knudsen moduli space $\overline{M}_{0,n}$; Thaddeus (cf. \cite{Thaddeus}) has realized the space of complete collineations as Chow quotients of suitable Grassmannians by suitable $\C^*$-actions. 

More recently, the authors, together with Occhetta, have studied the Chow quotients of flag varieties by the action of the maximal torus of their automorphism group, obtaining a full description in the case of the flag variety parametrizing complete flags of linear subspaces of $\C^4$ (cf. \cite{BOS2}), the variety of isotropic flags in $\C^4$ with respect to a symplectic form (cf. \cite{BiS}), and certain rational homogeneous spaces of Picard rank $1$ (cf. \cite{OS26}).  Although the approach introduced in \cite{BOS2} can be applied to a large class of varieties with torus actions (cf. \cite{OS26}), the examples studied show that a complete description of a Chow quotient and of its geometric properties can be extremely involved. In fact, Chow quotients dominate every GIT quotient, and their Picard number may become very large. Moreover, they are often singular and need not be Mori dream spaces. Consequently, a complete description of their geometry can be extremely difficult, and in the setting considered here such descriptions are known only in the above-mentioned  cases. 

The motivation behind this work was to find geometric quotients that retain enough information to reconstruct all the GIT quotients of a torus action, without having to describe the Chow quotient explicitly. This leads us to introduce the notion of a \emph{maximal geometric quotient} (\emph{MG quotient}, for short; cf. Definition \ref{def:MGquot}): a geometric quotient of $X$ such that every GIT quotient of $X$ is obtained from it by a rational contraction. In this way, the birational geometry of the whole collection of GIT quotients can be studied through the birational geometry of a single variety. MG quotients do not exist in general (cf. Remark~\ref{remark:MGnotExistence}), and, when they exist, they are unique only up to small $\Q$-factorial modifications.

This point of view is particularly useful when an MG quotient is a Mori dream space. Indeed, its Mori chamber decomposition encodes its rational contractions, and we may consider its {\em envelope}, defined as the normalization of the graph of its small $\Q$-factorial modifications (cf. Section ~\ref{ssec:envelopes}). We show that there exists a birational morphism from the envelope to the normalized Chow quotient (cf. Proposition~\ref{proposition:EnvChow}). Thus, the envelope provides a natural link between the birational geometry of an MG quotient and the geometry of the Chow quotient.

Our main result shows that, for the action of a maximal torus on a complete flag variety, this approach is particularly effective.

\begin{theorem}\label{theorem:Introduction}
	Let $X=G/B$ be a complete flag variety of dimension $n\geq 4$, with $G$ simple. Consider the action of the maximal torus $H\subset G$ by left multiplication. Then an MG quotient $\MX$ of $X$ by $H$ exists, and it is a Mori dream space with Picard rank equal to $2\rk(H)$. Moreover,  there exists an isomorphism of $\C$-algebras $\Cox(X)\simeq \Cox(\MX)$, and the birational map $\Env(\MX)\to \CX$ is an isomorphism.
\end{theorem} 

Another important feature of Theorem \ref{theorem:Introduction} is the relation between the Cox rings of $X$ and of its MG quotient. Although their Picard numbers are different, we prove that their Cox rings are isomorphic as complex algebras. The key geometric ingredient is the existence of a GIT quotient whose semistable locus has complement of codimension at least two, and on which the torus acts freely. This allows us to identify the divisor class group of the quotient with the group of equivariant line bundles on $X$, and consequently to identify their Cox rings, taking into account the natural equivariant grading. Furthermore, the last claim of the statement tells us that the rational contractions of the MG quotient are all the GIT quotients of $X$. The proof of the above Theorem will be the content of Section \ref{sec:mainproof}.

Finally, our approach also provides an explicit construction of an MG quotient of $X$, which we address in Section \ref{sec:ConsMG}. We start with a {\em Bruhat open neighborhood} of $X$ (see Section \ref{ssec:CoxRH}), whose MG quotients can be described as toric varieties of Picard number equal to the rank of the acting torus. We then show that an appropriate choice of such a quotient allows us to obtain an MG quotient of $X$ by performing exactly as many divisorial extractions as the rank of the torus, corresponding to the fundamental divisors contained in the complement of the Bruhat neighborhood. These extractions can be described as blowups along divisorial valuations and, remarkably, they can be performed in any order. This provides a concrete construction of the MG quotients, relating their birational geometry to the combinatorics of the root system and to the GIT chamber decomposition. 

\section{Preliminaries}\label{sec:prelim}

\subsection{Notation}\label{ssec:notation}

We work over the field of complex numbers. Let $M$ be a normal projective variety. We denote by $\NU(M)$ the $\R$-vector space of Cartier divisors on $M$ modulo numerical equivalence, and define the {\em Picard number} of $M$ as $\rho_M:=\dim \NU(M)$. A {\em contraction} of $M$ is a surjective morphism $\phi\colon M\to M'$ with connected fibers, where $M'$ is a normal projective variety. A birational map $\phi:M\dashrightarrow M'$, with $M'$ normal projective, is a {\em small modification} if it provides an isomorphism among two open sets whose complements have codimension at least two in $M$ and $M'$.  If $M$ and $M'$ are assumed to be $\Q$-factorial, such a map is called a {\em small $\Q$-factorial modification} (SQM, for short).  By a {\em rational contraction} of a $\Q$-factorial projective variety $M$ we mean a rational map $\phi\colon M \dashrightarrow M'$, with $M'$ normal projective, which factorizes as a composition of small $\Q$-factorial modifications followed by a contraction. Given a quasi-projective variety $X$ and an ample line bundle $L$ on $X$, we call $(X,L)$ a \emph{polarized pair}.  Finally, given a field $k$ and a lattice $L$, we denote by $L_k:=L\otimes_\Z k$ its associated $k$-vector space and, given a lattice homomorphism $f$, we denote by $f_k$ its corresponding $k$-linear map.

\begin{definition}\label{def:MDS}
A normal $\Q$-factorial projective variety $M$ is called a {\it Mori dream space} (MDS, for short) if \begin{enumerate}[leftmargin=\maren pt]
	\item\label{item:1} the Picard group $\Pic(M)$ of $M$ is finitely generated;
	\item\label{item:2} the nef cone $\Nef(M)$ is generated by classes of finitely many semiample divisors;
	\item\label{item:3} there exists a finite number of SQMs $f_i:M\dashrightarrow M_i$, for $i=0,\ldots,k$ with $M_0=M, f_0=\id_M$, such that every $M_i$ satisfies $(\ref{item:2})$ and $$\Mov(M)=\bigcup_{i=0}^k f_i^*\Nef(M_i).$$
\end{enumerate}
\end{definition}

\begin{definition}\label{def:Cox}
Let $M$ be a normal $\Q$-factorial projective variety with finitely generated class group $\Cl(M)$ and Picard number $\rho=\rho_M$. We will only deal with varieties with free class  group, so we assume $\Cl(M)\simeq\Z^\rho$ and choose Weil divisors $D_1,\ldots,D_\rho$ whose classes form a $\Z$-basis of $\Cl(M)$.
Then the {\em Cox ring of $M$} (with respect to $D_1,\ldots,D_{\rho}$) is the graded ring
\[ \Cox(M) := \bigoplus_{(m_1,\ldots,m_\rho)\in \Z^\rho} {\textstyle\HH^0\big(M,\cO_M\big(\sum_{i=1}^\rho m_iD_i\big)\big)}.\]
The ring structure of $\Cox(M)$ is induced by the multiplication of sections in $\C(M)$. The {\em characteristic torus of $M$} is defined as $\Hom(\Cl(M), \C^*)=\Spec(\C[\Cl(M)])$, so its character lattice $\Cl(M)$ is the grading group of $\Cox(M)$.
\end{definition}

\begin{remark}
Note that the isomorphism class of the $\C$-algebra $\Cox(M)$ is independent of the chosen basis. Moreover, by \cite[Proposition~2.9]{HuKeel}, we know that a variety $M$ as in Definition \ref{def:Cox} is an MDS if and only if $\Cox(M)$ is a finitely generated $\C$-algebra.
\end{remark}

\subsection{Envelopes}\label{ssec:envelopes}

We introduce the notion of envelope of an MDS, which can be intuitively understood as a normal projective variety which resolves all the SQMs of the MDS. We remark that a similar notion was
first discussed in \cite[Remark~4.6]{BakerChow} with the name \textit{Mori master space}, and it was also considered in \cite[\S~3]{PU} as a universal model for an MDS.

\begin{definition}\label{definition:Envelope}
	Let $M$ be an MDS, and  consider its SQMs $f_i:M\dashrightarrow M_i$, for $i=0,\ldots,k$, as in Definition \ref{def:MDS}. 
	Let $U$ be a nonempty open subset of $M$ contained in the locus where all the $f_i$'s are isomorphisms; we consider the graph
	\[ \Gamma_{M} = \overline{\{(f_0(x),f_1(x),\ldots,f_k(x)) \mid x\in U\}} \subset \prod_{i=0}^k M_i.\]
	The \textit{envelope} $\Env(M)$ is defined as the normalization of $\Gamma_{M}$.
\end{definition}

\begin{remark}\label{remark:Envelope}
	We can also describe the envelope as the projective spectrum of a finitely generated $\C$-algebra. For each $i$, let $A_i$ be a sufficiently ample
Cartier divisor on $M_i$, and denote by
$D_i:=f_i^*A_i$ its strict transform on $M$. For any $m\geq 0$,  consider 
	\[ \Phi_m  \colon \bigotimes_{i=0}^k \HH^0(M,\cO_M(mD_i)) \lra \HH^0(M,\cO_M(m(D_0+\ldots+D_k))),\]
	and let $R_m$ be the image of $\Phi_m$. Then one can show that $R:= \bigoplus_{m\geq 0} R_m$ is a finitely generated $\C$-algebra, and $\Gamma_M = \Proj R$. By setting $\overline{R}$ to be the integral closure of $R$ in the section ring $R(M; \cO_M(D_0+\ldots+D_k))$, we obtain that 
	$$\Env(M)=\Proj \overline{R}.$$
\end{remark}

\begin{remark}
	In the language of torus actions and their quotients, that we will introduce in the following section, we remark that the envelope of an MDS $M$ can be also described as the limit quotient of $\Spec(\Cox(M))$ by the action of the characteristic torus of $M$.
\end{remark}

\subsection{Torus actions and  GIT quotients}\label{ssec:torusquots}

Let $X$ be a normal quasi-projective variety, endowed with  a non-trivial action of a complex torus $H$ of rank $r$. We denote by  $\Mo(H)$ the character lattice of $H$. 
The fixed point locus decomposes as $X^H=\bigsqcup_{Y\in \cY} Y$, where each $Y$ is a connected component of $X^H$. Note that (cf. \cite[Theorem~1.1]{IVERSEN}), whenever $X$ is smooth, each $Y$ is smooth, hence irreducible.

Given a line bundle $L$ on $X$, a {\em linearization} of $L$ by $H$ is an $H$-action $H\times L\to L$ which is linear along the fibers and such that the projection $L\to X$ is $H$-equivariant; we will say that the linearization defines an action of $H$ on the pair $(X,L)$. Note that linearizations always exist \cite[Remark~2.4]{KKLV}, and two linearizations of $L$ differ by a character of $H$.

For every $p\in Y$ of a fixed point component $Y\in\cY$, we denote by $\mu_L(Y)\in \Mo(H)$ the {\em weight} of the $H$-action on the fiber $L_p$, which depends only on the connected component $Y$. We then denote by $\Delta$ the {\em polytope of fixed points of the action}, that is the convex hull in $\Mo(H)_{\R}$ of $\mu_L(Y)$, for $Y\in \cY$. The action on  $(X,L)$ induces a linear action of $H$ on $\HH^0(X,L)$, and a subsequent decomposition of this space in eigenspaces associated to weights. We denote by $\Gamma$ the {\em polyhedron of sections}, 
defined as the convex hull in $\Mo(H)_{\R}$ of the weights $\tau\in \Mo(H)$ such that $\HH^0(X,L)_{\tau}\neq 0$. If $X$ is projective, then $\Gamma$ is bounded and we will refer to it as the polytope of sections of $L$.

Given an ample line bundle $L$ on $X$ and $u\in  \Gamma$, we set
$$\cA_u := \bigoplus_{m\geq 0, mu\in \Z^r}\HH^0(X,mL)_{m u}.$$
In the projective and affine settings that we will consider here, the $\C$-algebra $\cA_u$ is finitely generated, and we define the \emph{GIT quotient} of $(X,L)$ under the $H$-action as the normal quasi-projective variety
\[
\GX_u:=\Proj\cA_u.
\]
This is a good quotient of a semistable open subset $X^{\ss}_u\subset X$, determined by the linearization associated with the weight $u$. We call {\em GIT quotient of $X$} any GIT quotient obtained from a polarized pair $(X,L)$.

\begin{example}[MDS's and VGIT]\label{ex:VGIT}
A key example of the above situation $X$ is an MDS, then one can consider the characteristic torus action on $\Spec(\Cox(X))$. It was shown by Hu and Keel \cite{HuKeel} that there is a one-to-one correspondence among rational contractions of $X$  and GIT quotients of the characteristic torus action on the affine variety $\Spec(\Cox(X))$ linearized on the trivial line bundle.
\end{example}

Finally, we fix some notation in the case of torus actions of rank $1$, and refer to \cite{WORS1} for details. Let $\C^*$ act nontrivially on $X$.  For any point $x\in X$, we denote by $\lim_{t\to 0}tx,\lim_{t\to \infty} tx$ the limit points of its orbit closure. Moreover, for any fixed connected component $Y\subset X^{\C^*}$ we define the {\it Bia\l ynicki-Birula cells} (shortly, BB-cells) as
$$ X^+(Y)=\{x\in X\mid \lim_{t\to 0} tx\in Y\}, \quad X^-(Y)=\{x\in X\mid \lim_{t\to \infty} tx\in Y\}.$$
If $X$ is smooth and projective, by the Bia\l ynicki-Birula theorem (cf. \cite{BB}, see also \cite{BBC}), there exists a unique fixed connected component $Y_-$ (resp. $Y_+$) such that $X^-(Y_-)$ (resp. $X^+(Y_+)$) is dense; we call such component $Y_-$ (resp. $Y_+$) the {\em sink} (resp. the {\em source}) of the $\C^*$-action on $X$. 
We refer to $Y_{\pm}$ as the {\em extremal fixed point components} of the action; the other connected fixed point components will be called {\em inner}. 

\subsection{Chow and limit quotients} 

We briefly recall the notion of Chow quotient and limit quotient, following respectively \cite{Kapranov} and \cite{BakerChow}.

\subsubsection{Chow quotients} Let $X\subset \P^n$ be a smooth projective variety, endowed with  a non-trivial action of an algebraic torus $H$. Up to shrinking, one may find an open $H$-invariant subset $U$ of $X$ such that, for any $z\in U$, the orbit closures $\overline{H\cdot z}$ have constant dimension $k$ and degree $d$. Then, there exists an embedding
\[ \psi: U/H\hookrightarrow \text{Chow}_{k,d}(X),\quad [H\cdot z]\mapsto [\overline{H\cdot z}]\]
from the geometric quotient $U/H$ to the Chow variety of $X$ parametrizing cycles of $X$ of dimension $k$ and degree $d$.

\begin{definition}
	The \textit{normalized Chow quotient} $\CX$ is the normalization of the closure of $\psi(U/H)$ in $\text{Chow}_{k,d}(X)$.
\end{definition}

If no confusion arises in the sequel of the paper we will commonly refer to $\CX$ as the Chow quotient of $X$ by $H$. Recall that the Chow quotient dominates every GIT quotient of $X$.

\begin{lemma}\cite[Theorem~(0.4.3)]{Kapranov}
	For any GIT quotient $\GX$ of $X$ by $H$ there exists a birational morphism $\CX\to \GX$.
\end{lemma}

\subsubsection{Limit quotients}

Let $X$ be a normal projective variety endowed with  a non-trivial  action of a torus $H$. The set of GIT quotients of $X$ forms a finite inverse system. There exists a unique component in the inverse limit of this system dominating every GIT quotient of the system. The normalization of this component will be called the {\em limit quotient} and denoted  by $\LX$. In other words, for every GIT quotient $\GX$ we have a contraction \[\LX\longrightarrow\GX\] and $\LX$ is categorically final for this property.

It is known (see \cite[Corollaries~2.6,~2.7]{BakerChow}) that $\LX$ can be computed out of the subsystem of GIT quotients of the polarized pair $(X,L)$, for any ample line bundle $L$, and that, in the case in which $X$ is projective, it coincides with the Chow quotient of $X$.

\subsection{MG quotients}\label{ssec:MG}

Let us go back to the case of a normal complex projective variety $X$ endowed with the action of a complex torus $H$, and the corresponding inverse system of GIT quotients of $X$. Different concepts of quotient lying in between Chow and GIT can be already found in the literature (see for instance \cite[Definition~3.12]{HHW}, \cite{Gh24}); our notion is motivated by the well-known  relation among birational geometry and variation of GIT. More precisely, given a GIT quotient $\GX$, the associated collection of rational contractions of $\GX$ contains some other GIT quotients of $X$. It then makes sense to try to find a GIT quotient whose system of rational contractions contains as much information of the GIT system as possible. 

\begin{definition}\label{def:MGquot}
Let $X$ be a normal quasi-projective  variety endowed with the action of a torus $H$. Let $L$ be an ample line bundle on $X$. A {\em maximal geometric quotient} ({\em  MG quotient}, for short) of the polarized pair $(X,L)$ (resp. of $X$)  is a $\Q$-factorial geometric quotient $\MX$ of $(X,L)$ (resp. of $X$) such that  for every GIT quotient $\GX$ of $(X,L)$ (resp. of $X$) the natural rational map \[\MX\dashrightarrow \GX\] is a rational contraction.
\end{definition}

In the cases which  MG quotients exist, they are related to the corresponding Chow quotients via the envelope construction introduced above. In fact we may state the following.

\begin{proposition}\label{proposition:EnvChow}
	Let $(X,L)$ be a polarized pair, with $X$ a normal projective variety, endowed with the action of a torus $H$. Assume that there exists an  MG quotient $\MX$ of the pair $(X,L)$ that is an MDS. Then  there exists a birational morphism \[\Env(\MX)\longrightarrow\CX.\]
\end{proposition}

\begin{proof}
It is enough to note that under the above assumptions, the inverse system of rational contractions of $\MX$ contains the system of GIT quotients of $(X,L)$, and so the claim follows.
\end{proof}

In the above situation, the morphism $\Env(\MX)\to\CX$ is not necessarily an isomorphism. In other words, in general not all the rational contractions of $\MX$ correspond to quotients of $X$.

\begin{remark}\label{remark:MGnotExistence}
 MG quotients of torus actions do not exist, in general. An example can be easily constructed in the setup of equalized $\C^*$-actions on smooth varieties, as presented in \cite{WORS5}. Consider the case of a polarized pair $(X',L')$, with $X'$ a smooth projective variety endowed with an equalized $\C^*$-action. Assume that the criticality $r$ of the action is at least three, and that the ranks of the positive and negative parts of the normal bundles of every inner fixed point component of the action  are larger than one (examples of this kind can be found in \cite{WORS5, FS}). 
Under these hypotheses we have $r$ geometric quotients of $(X',L')$, that are smooth projective varieties, related by Atiyah flips.
\[
\GX'(0,1)\dashrightarrow \GX'(1,2)\dashrightarrow \dots \dashrightarrow \GX'(r-2,r-1) \dashrightarrow \GX'(r-1,r)
\]

Let $Y_1,Y_2$ be two inner fixed point components of the action of different weights, corresponding to indices $i_1,i_2\in\{1,\dots,r-1\}$. Let $\pi:X\to X'$ be the blowup of $X'$ along $Y_1$ and $Y_2$, endowed with an ample line bundle of the form $L=\pi^*(mL')-a_1E_1-a_2E_2$, where $a_1,a_2>0$ and $m\gg 0$. Then $(X,L)$ has $r+2$ geometric quotients, $r$ of them equal to the geometric quotients of $X'$; the two extra ones, that we denote by $\GX_1$ and $\GX_2$, dominate two geometric quotients of $X'$ via smooth blowups:
\[\begin{tikzcd}
                                     & \GX_j \arrow[ld, "\mathrm{blowup}"'] \arrow[rd, "\mathrm{blowup}"] &                   \\
{\GX'(i_j-1,i_j)} \arrow[rr, dashed] &                                                                    & {\GX'(i_j,i_j+1)}
\end{tikzcd}\]
Consequently, no geometric quotient of $X$ admits rational contractions to both $\GX_1$ and $\GX_2$.
\end{remark}  

\section{MG quotients of flag varieties}\label{sec:MGRH}

In this section, we discuss in detail the concepts introduced in Section \ref{sec:prelim} in the special case of rational homogeneous varieties and prepare the ground for the proof of the main theorem, discussed in Section \ref{sec:mainproof}. More precisely, we will work in the following setup.

\begin{setup}\label{set:flags}
Let  $G$ be a connected algebraic group, $B\subset G$ a Borel subgroup, $H\subset B$ a maximal torus and $X=G/B$ the corresponding complete flag variety, endowed with the action of the torus $H$. We will assume that the Lie algebra $\fg$ of $G$ is simple and that the Dynkin diagram of $G$ is neither $\DA_1$ nor $\DA_2$. 
\end{setup}

The two excluded cases are discussed in Remark \ref{rem:main}.

\subsection{The Cox ring of a complete flag variety}\label{ssec:CoxRH}

Under the hypotheses of Setup~\ref{set:flags}, we may assume, without loss of generality, that $G$ is the {\em adjoint group} of $\fg$, so that the {\em character lattice} $\Mo(H)$ of $H$ equals the {\em root lattice} of  $\fg$. The Borel subgroup $B\subset G$ determines a set of {\em positive roots} $\Phi^+$, with their corresponding {\em simple roots} $\{\alpha_1,\dots,\alpha_n\}$ generating the lattice of $\Mo(H)$, where $n=\rk(H)$.

Let us consider also a maximal torus $\tilde H$ of the simply connected group of $\fg$, whose lattice of characters $\Mo(\tilde H)$ equals the {\em weight lattice} of $\fg$. Moreover, there is a well-known isomorphism of groups
\[\Mo(\tilde H)\simeq\Pic(X),\] in which a weight $\omega \in \Mo(\tilde H)$ corresponds to a line bundle $L(\omega)$ on $X$. In this way, we can define the $\C$-vector space on which the Cox ring of $X$ is supported as
\[
\Cox(X):=\bigoplus_{\omega\in \Mo(\tilde H)}\HH^0\big(X,L(\omega)\big)=\bigoplus_{a_1,\dots,a_n\in\Z_{\geq 0}}{\textstyle\HH^0\big(X,L\big(\sum_{i=1}^n a_i\omega_i\big)\big)}.
\]
The natural $M(\tilde H)$-graded $\C$-algebra structure of $\Cox(X)$ comes from the  product of the function field $\C(X)$.  Here $\omega_i\in \Mo(\tilde H)$, $i\in\{1,\dots,n\}$ denote the {\em fundamental weights} of $\fg$, defined by the property $2(\omega_i,\alpha_j)/(\alpha_j,\alpha_j)=\delta_{ij}$, where $(\cdot,\cdot)$ denotes the pairing in $\Mo(H)_\R=\Mo(\tilde H)_\R$ induced by the {\em Killing form} of $\fg$. In other words, the line bundles $L(\omega_i)$ are the pullbacks of the ample generators of the Picard groups of the $n$ contractions of $X$ of Picard number one. 
We will denote the {\em Cartan matrix} of $G$ by 
\[C:=\left(\langle \alpha_i,\alpha_j\rangle\right)_{i,j},\quad\textrm{where }\langle \alpha,\beta\rangle:=2\frac{(\alpha,\beta)}{(\beta,\beta)}\quad\textrm{for all }\alpha,\beta\ne0.\]

The fixed point locus of $X$ by the $H$-action is the set of points $\{\sigma B\mid\sigma\in W\}$, where $W:=N_G(H)/H$ is the {\em Weyl group} of $G$ with respect to $H$. The Weyl group is usually identified with the group generated by the {\em simple reflections} associated to $\alpha_1,\dots,\alpha_n$, namely \[r_i:\Mo(H)_\R\to\Mo(H)_\R,\quad r_i(\alpha):=\alpha-\langle\alpha,\alpha_i\rangle\alpha_i.\] 
Every $H$-fixed point $\sigma B$ admits an $H$-invariant open neighborhood, $H$-equivariantly isomorphic to $T_{X,\sigma B}$, 
that we call the {\em Bruhat open neighborhood $U_\sigma\subset X$ of $\sigma B\in X$}. It can be defined as the isomorphic image of the nilpotent Lie subalgebra $\Ad_\sigma (\fn)\subset \fg$, where $\fn:=\bigoplus_{\alpha\in \Phi^+}\fg_{-\alpha}$, by the $H$-equivariant map
\[x\in \Ad_\sigma (\fn)\mapsto \sigma\exp(\Ad_{\sigma^{-1}}(x))B\in X.\]

Following \cite[Corollary~3.8]{WORS5}, every line bundle $L(\omega)$ admits an $H$-linearization whose weights at the $H$-fixed points are
\begin{equation}\label{eq:weights}
\mu_{L(\omega)}(\sigma B)=\omega-\sigma(\omega)\in\Mo(H), \mbox{ for every }\sigma\in W.
\end{equation}
This defines a splitting of the group of {\em $H$-linearized line bundles} in $X$, that is
\begin{equation*}
\Pic^H(X)\simeq	\Pic(X)\oplus \Mo(H)\simeq \Mo(\tilde H)\oplus \Mo(H)\simeq \Mo(\tilde H\times H).
\end{equation*} 
In the sequel, we will always identify $\Pic^H(X)$ with $\Mo(\tilde H\times H)$.

With this choice of $H$-linearizations, if $s\in  \HH^0(X,L(\omega))_m, s'\in  \HH^0(X,L(\omega'))_{m'}$, we have that $ss'\in \HH^0(X,L(\omega+\omega'))_{m+m'}$. This gives a grading of $\Cox(X)$ over the character lattice of the torus $\tilde H\times H$, namely 
\[
\Cox(X)=\bigoplus_{(\omega,m)\in \Pic^H(X)}\HH^0\big(X,L(\omega)\big)_m
\]
hence an action of the torus $\tilde H\times H$ on the affine variety $\Spec(\Cox(X))$.

\begin{remark}\label{rem:equivcox}
Note that, as a $\Pic^H(X)$-graded algebra, $\Cox(X)$ can be seen as a graded subalgebra of the {\em equivariant Cox ring} $\Cox^H(X)$ introduced in \cite{Vez} (see also \cite{Br07}), namely its $H$-invariant part. Instead, our situation provides an {\em equivariant Cox construction} of $X$, in the sense of \cite[\S 4]{ArzHau}.
\end{remark}

\begin{definition}\label{def:effconeofMG}
The {\em weight cone} of $X$ is the convex cone defined as \[\Delta:=\R_{\geq0}(\{(\omega,m)\in\Pic^H(X)_\R\mid \HH^0(X,L(\omega))_m\neq0\}).\]
\end{definition}

\begin{remark}\label{rem:coxRH}
Given a weight $\omega=\sum_ia_i\omega_i$, $a_i\geq 0$, the line bundle $L(\omega)$ is globally generated; therefore, by \cite[Lemma~2.4~(3)]{BWW} and considering the linearization of $L(\omega)$ introduced above, we have that
\[
\Delta(\omega):=\{m\in \Mo(H)_\R\mid (\omega,m)\in \Delta\}
\]
equals both the weight polytope and the polytope of sections of the line bundle $L(\omega)$. Note that $\Delta(\omega)$ can be written as a Minkowski combination of the  polytopes of sections of the fundamental weights, that is 
\[\textstyle \Delta(\omega)=\Delta(\sum_ia_i\omega_i)=\sum_ia_i\Delta(\omega_i).\]
Note also that, by definition, the cone
\[\{\omega\in \Mo(\tilde{H})_\R\mid \Delta\cap \big(\{\omega\}\times \Mo(H)\big)\neq \emptyset\}\]
equals both the nef and the effective cone of $X$.
\end{remark}

In particular, with the above notation, 
\[ \Delta=\R_{\geq 0}\{(\omega_i,\omega_i-\sigma(\omega_i))\mid i\in\{1,\dots,n\},\,\, \sigma\in W\}. \]

\begin{remark}\label{rem:polyomegai}
Let us observe that the general theory of representations of semisimple groups tells us that the vertices of the polytope $\Delta(\omega_i)$ are the elements in the $W$-orbit of $-\omega_i$, translated by $\omega_i$. Since the isotropy subgroup of $\omega_i$ by the $W$-action is the subgroup $W^i\subset W$ generated by the reflections $\{r_j\mid j\neq i\}$, we may identify the set of vertices of $\Delta(\omega_i)$ with the set $W/W^i$ of left cosets modulo $W^i$. This can be interpreted as follows: the linear system of global sections of $L(\omega_i)$ maps the flag variety $G/B$ to a generalized Grassmannian, i.e., a rational homogeneous variety $G/P^i$ of Picard number one. Then $\Delta(\omega_i)$ coincides with the weight polytope of $G/P^i$, polarized with the line bundle whose pullback to $G/B$ is $L(\omega_i)$, and the vertices of this polytope are the weights of the action at the points $\sigma P^i$, $\sigma\in W$, namely $\omega_i-\sigma(\omega_i)$. It is then enough to note that $\sigma P^i=\sigma'P^i$ if and only if $\sigma^{-1}\sigma'\in W^i$. Finally, by looking at the action of $H$ on the tangent space $T_{G/P^i,\sigma P^i}$, we obtain a local description of $\Delta(\omega_i)$ around $\omega_i-\sigma(\omega_i)$: a neighborhood of this point in $\Delta(\omega_i)$ equals a neighborhood of the cone with vertex $\omega_i-\sigma(\omega_i)$ passing by the points $\omega_i-\sigma(\omega_i-\beta)$, $\beta\in \Phi^+(i)$, where $\Phi^+(i)$ is the set of positive roots having a positive $i$-th coordinate (with respect to the base of simple roots). 
\end{remark}

\subsection{Fundamental $\C^*$-actions}\label{ssec:fund}

In this section we will  geometrically identify the divisors corresponding to elements in the space $\HH^0(X,L(\omega_i))_{\omega_i-\sigma(\omega_i)}$, $\sigma\in W$. To this end, we will make use of the following actions.
\begin{definition}\label{def:fundaction}
For every $i\in\{1,\dots,n\}$, $\sigma\in W$, we consider the {\em fundamental $\C^*$-action} determined by the cocharacter
\[
\mu_i^\sigma:\Mo(H)\to \Z,\quad \mu_i^\sigma(\sigma^{-1}(\alpha_j)):=\delta_{ij}.
\] 
For simplicity, we will write $\mu_i:=\mu_i^e$, so that $\mu_i^\sigma=\mu_i\circ\sigma$ for every $i,\sigma$. We abuse notation by denoting also by $\mu_i^\sigma$ both their associated $\R$-linear maps and $\C^*$-actions. \end{definition}

\begin{remark}\label{rem:fundcochar}
Observe that we may write every $\omega\in \Mo(H)_\R$ as 
\begin{equation*}
\label{eq:coordinates} \omega=\sum_{i=1}^n\langle \omega,\alpha_i\rangle \omega_i=\sum_{i=1}^n\mu_i(\omega)\alpha_i.
\end{equation*}
In particular, for every $j\in\{1,\dots,n\}$ we have the relations
\begin{equation*}
\label{eq:coordinates2}
\langle \omega,\alpha_j\rangle =\sum_{i=1}^n \mu_i(\omega)\langle \alpha_i,\alpha_j \rangle.
\end{equation*}
\end{remark}

Note now that, with our prescribed linearization (\ref{eq:weights}), for every $i,\sigma$ the fundamental $\C^*$-action $\mu_i^\sigma$ acts on $(X,L(\omega_i))$ with weights of the form $\mu_i^\sigma(\omega_i-\tau(\omega_i))$, compatibly with the immersion $G/P^i\subset \P(\HH^0(X,L(\omega_i))^\vee)$. Furthermore, we may state the following result.

\begin{lemma}\label{lem:omegaifundi}
Consider the fundamental $\C^*$-action $\mu_i^\sigma$ on $X$. There exists a unique fixed point component $Y_i^\sigma$ satisfying that $\ol{X^-(Y _i^\sigma)}$ is a prime divisor. Its weight with respect to the line bundle $L(\omega_i)$ is $\mu_i^\sigma\mu_{L(\omega_i)}(Y^\sigma_i)=\mu_i^\sigma(\omega_i)-\mu_i(\omega_i)+1$.
\end{lemma}

\begin{proof}
Let us consider the $\C^*$-action $\mu_i^\sigma$ induced on the generalized Grassmannian $G/P^i\subset\P(\HH^0(X,L(\omega_i))^\vee)$. The fixed point components of this action are the images of the fixed point components of the action on $X$, and the same property holds for the corresponding negative BB-cells. Since the morphism $X\to G/P^i$ is smooth, it follows that the closure of a BB-cell in $X$ is a divisor if and only if its image in $G/P^i$ is a divisor. 

From the description of $\Delta(\omega_i)$ presented in Remark \ref{rem:polyomegai}, we get that any possible vertex of the polytope is of the form $\omega_i-\sigma^{-1}(\omega_i-\sum_j a_j\alpha_j)$, with $a_i\geq 0$ necessarily. Then its projection via $\mu_i^\sigma$ has a minimum $\mu_i^\sigma(\omega_i-\sigma^{-1}(\omega_i))$, which is achieved only at the point $\sigma^{-1} P^i$.

Note also that $\Delta(\omega_i)$ contains the point $\omega_i-\sigma^{-1}(r_i(\omega_i))=\omega_i-\sigma^{-1}(\omega_i-\alpha_i)$, whose projection via $\mu_i^\sigma$ is $\mu_i^\sigma(\omega_i-\sigma^{-1}(\omega_i))+1$. This shows that we have a fixed point component $Z_i^\sigma$ in $G/P^i$ of that precise weight. Now, since the sink of the action is isolated, we necessarily have, by Bend-and-Break, that
the rank of the positive summand of the normal bundle of $Z_i^\sigma$ in $G/P^i$ is one;  equivalently, $\ol{(G/P^i)^-(Z_i^\sigma)}$ is a divisor in $G/P^i$ (cf. \cite[Lemma~4.5]{WORS3}). By \cite[Lemma~2.9~(2)]{WORS1}, since the Picard number of $G/P^i$ is one, an inner fixed point component of this kind is unique. Finally, by our previous argument, there exists a unique irreducible component $Y_i^\sigma$ in $X$ dominating $Z_i^\sigma$ satisfying the required properties. 
\end{proof}

\begin{remark}\label{rem:omegaifundipt}
From the arguments in the proof, we may also describe $Y_i^\sigma$ as the unique $\mu_i^\sigma$-fixed point component containing the point $\sigma^{-1}r_iB$. In particular, for $j\neq i$, the $L(\omega_j)$-weight of $Y_i^\sigma$ is equal to 
\[\mu_i^\sigma\mu_{L(\omega_j)}(Y^\sigma_i)=\mu_i^\sigma\big(\omega_j-\sigma^{-1}r_i(\omega_j)\big)=\mu_i^\sigma\big(\omega_j-\sigma^{-1}(\omega_j)\big)=\mu_i^\sigma(\omega_j)-\mu_i(\omega_j).\]
In particular, given $\omega=\sum_{j=1}^n\langle \omega,\alpha_j\rangle\omega_j\in \Mo(\tilde H)$, we get that the $L(\omega)$-weight of the $\mu_i^\sigma$-fixed point component $Y_i^\sigma$ is equal to
\[\begin{array}{rl}\mu_i^\sigma\mu_{L(\omega)}(Y_i^\sigma)&=\mu_i^\sigma\mu_{L(\omega)}(\sigma^{-1}r_iB)=\\[4pt]
&=\sum_{j=1}^n \langle\omega,\alpha_j\rangle\big(\mu_i^\sigma(\omega_j)-\mu_i(\omega_j)+\delta_{ij}\big)=\\[4pt]
&=\mu_i^\sigma(\omega)-\mu_i(\omega)+\langle\omega,\alpha_i\rangle=\mu_i^\sigma(\omega)-\mu_i^{r_i}(\omega).\end{array}
\]
\end{remark}

\begin{remark}\label{rem:omegaifundi}
Taking into account that, denoting by $w_0\in W$ the longest element of the Weyl group, there exists a permutation $p$ of $\{1,\dots,n\}$ (depending on the Dynkin diagram of $G$) such that
\[\mu_i^{w_0 \sigma}=-\mu_{p(i)}^{\sigma},\quad\mbox{for every }\sigma\in W,\]
we may conclude that the sets of closures of positive and negative BB-cells of fundamental $\C^*$-actions are equal.
\end{remark}

\subsection{Fundamental divisors}\label{ssec:funddiv}

Our next step will be to show that the BB-cells of the $\C^*$-action on $X$ induced by any cocharacter $\mu:\Mo(H)\to\Z$ are intersections of BB-cells of fundamental $\C^*$-actions $\mu_i^\sigma$. In particular, we will be able to conclude that the only BB-cells of $\C^*$-actions in $X$ that have codimension one are the ones provided in Lemma \ref{lem:omegaifundi}.

\begin{remark}\label{rem:Weylfan}
Note that the normal fan $\Sigma_W\subset \No(H)_\R$ of the weight polytope $\Delta(\omega)$ associated to  an ample line bundle $L(\omega)$ in $X$ is independent of $\omega$, and we call it the {\it Weyl fan of $X$}. 
The lattice generators of the rays of $\Sigma_W$ are precisely the cocharacters $\mu_i^\sigma\in \No(H)$, and the cones of maximal dimension of $\Sigma_W$ are in one-to-one correspondence to vertices of $\Delta(\omega)$, hence to elements $\sigma\in W$ (see also Section \ref{ssec:GITquotflag}). In particular, for every $\mu\in \No(H)$ there exist $\sigma\in W$ (not unique, in general), indices $i_1,\dots,i_k\in\{1,\dots,n\}$, and positive integers $a_1,\dots,a_k\in \Z$ such that \[\mu=\sum_{r=1}^ka_r\mu^\sigma_{i_r}.\]
\end{remark} 

\begin{proposition}\label{prop:BBcellsfund}
Let $\mu\in\No(H)$ be a cocharacter of $H$, and denote also by $\mu$ the corresponding $\C^*$-action on $X$. Then the BB-cells of $\mu$ are intersections of BB-cells of fundamental $\C^*$-actions on $X$.
\end{proposition}
 
\begin{proof}
Note first that, arguing as in Remark \ref{rem:omegaifundi}, we may restrict ourselves to the case of negative BB-cells of $\C^*$-actions $\mu$. 

Given $\mu\in \No(H)$, we write it as a positive combination $\mu=\sum_{r=1}^ka_r\mu^\sigma_{i_r}$ of fundamental $\C^*$-actions, for some $\sigma\in W$, and we will show that the BB-cells of $\mu$ are intersections of BB-cells of the $\C^*$-actions $\mu^\sigma_{i_r}$. This fact can be shown locally, by proving that the BB-cells of the action of $\mu$ on a Bruhat neighborhood $U_\tau$ are intersections of BB-cells of the actions of $\mu^\sigma_{i_r}$ on $U_\tau$, for every $\tau$.

By using the action of $W$, it is then enough to show that the above claim holds in the case $\tau=e$, for every possible choice of $\sigma$, $i_1,\dots,i_k$, $a_1,\dots,a_k>0$. In order to do so, we use the $H$-equivariant isomorphism $U_e\simeq \fn=\bigoplus_{\alpha\in\Phi^+}\fg_{-\alpha}$ and note that the negative BB-cells of $\mu$ in $\fn$ are linear subspaces of the form
\[\fn^-(0)=\bigoplus_{\mu(-\alpha)\leq 0}\fg_{-\alpha}.\]

Then we may conclude by noting that, if $\mu=\sum_{r=1}^ka_r\mu^\sigma_{i_r}$, $a_r>0$, we have that, for every root $\alpha$, 
\[\mu(-\alpha)\leq 0 \iff \mu^\sigma_{i_r}(-\alpha)\leq 0\mbox{ for every }r.\]
One implication is obvious. For the other, the set $\{\mu_{i_r}^\sigma(-\alpha)\mid r=1,\dots, k\}$ is either contained in $\Z_{\geq 0}$ or in $\Z_{\leq 0}$, for every $\alpha$. Then if $\mu(-\alpha)\leq 0$, we may easily conclude that $\mu^\sigma_{i_r}(-\alpha)\leq 0$ for every $r$. 
\end{proof} 

\begin{corollary}\label{cor:BBcellsfund}
The only closures of BB-cells of $\C^*$-actions in $X$ of codimension one  are the divisors $\ol{X^-(Y_i^\sigma)}$ presented in Lemma \ref{lem:omegaifundi}. \qed
\end{corollary}

\begin{definition}\label{def:fundamental}
With the notation of Lemma \ref{lem:omegaifundi}, the divisor $\ol{X^-(Y_i^\sigma)}$ in the complete flag variety $X=G/B$ will be denoted by $D_i^\sigma$, and called the {\em fundamental divisor} associated to the $\C^*$-action $\mu_i^\sigma$. 
\end{definition}

\begin{lemma}\label{lem:funddivi}
There is an isomorphism $\cO_X(D_i^\sigma)\simeq L(\omega_i)$.
\end{lemma}

\begin{proof}
Let us consider $G/P^i$ as a closed subset in $\P(\HH^0(X,L(\omega_i))^\vee)$. The point $\sigma^{-1}P^i\in G/P^i$ is the projectivization of the eigenspace $\HH^0(X,L(\omega_i))^\vee_{\omega_i-\sigma^{-1}(\omega_i)}$, which is $1$-dimensional. Now we consider the intersection $D$ of $G/P^i$ with the hyperplane defined by 
\[
\P\left(\bigoplus_{m\neq \omega_i-\sigma^{-1}(\omega_i)}\HH^0(X,L(\omega_i))^\vee_m\right).
\]
The divisor $D$ is $H$-invariant and, by definition, it does not meet the (open) BB-cell $X^-(\sigma^{-1}P^i)$. Hence it is contained in the union of the closure of the negative BB-cells of the $\mu_i^\sigma$-fixed point components different from the sink. But, as in the proof of Lemma \ref{lem:omegaifundi}, the only closure of a negative BB-cell that has codimension one is the one whose inverse image into $X$ is $D_i^\sigma$; hence it follows that the inverse image of $D$ into $X$ is $D_i^\sigma$. This completes the proof.
\end{proof}

From the above proof, we may immediately state the following corollary, that tells us that the weight cone $\Delta$ is generated by the weights of the action on the fundamental divisors $D_i^\sigma$.

\begin{corollary}\label{cor:funddivi}
The $H$-action on $\P(\HH^0(X,L(\omega_i)))$ induced by our chosen linearization of $L(\omega_i)$ has weight $\omega_i-\sigma^{-1}(\omega_i)$ on the divisor $D_i^\sigma$, for every $\sigma\in W$. \qed
\end{corollary}

\subsection{Weight polytopes of $X$}\label{ssec:GITquotflag}

Given $\omega\in\Q_{>0}(\omega_1,\dots,\omega_n)\subset\Mo(H)_\Q$, so that the rational line bundle class $L(\omega)$  is $\Q$-ample, the different GIT quotients of $(X,L(\omega))$ are determined by the GIT chamber decomposition of the weight polytope $\Delta({{\omega}})$. Let us start by writing this polytope in terms of the cocharacters $\mu_i^\sigma$ defining the fundamental $\C^*$-actions; a straightforward computation leads to  
\begin{equation*}
\Delta({{\omega}})=\left\{m\in \Mo(H)_\R\mid \mu_i^\sigma(m)\geq \mu_i^\sigma({{\omega}})-\mu_i({{\omega}}),\mbox{ for every }i,\sigma\right\}.
\end{equation*}

We now introduce the following relevant subset of $\Delta(\omega)$.

\begin{definition}\label{def:Deltamov}
We call {\em movable weight polytope} of $(X,L(\omega))$ the set
\[
\Delta^{\mov}({{\omega}}):=\left\{m\in \Mo(H)_\R\mid \mu_i^\sigma(m)\geq\mu_i^\sigma({{\omega}})-\mu_i^{r_i}({{\omega}}),\mbox{ for every }i,\sigma\right\}.
\]
We say that $\omega\in\Q_{>0}(\omega_1,\dots,\omega_n)$ is a {\em good weight} if the movable weight polytope $\Delta^{\mov}(\omega)$ has non-empty interior. 
\end{definition}

As we will see in the next section, good weights exist and provide MG quotients for every $X$.

\begin{remark}\label{rem:Deltamov} The inclusion $\Delta^{\mov}(\omega)\subseteq\Delta(\omega)$ follows from the fact that, for every $i$ we have $\mu_i^{r_i}(\omega)=\mu_i(\omega)-\langle\omega,\alpha_i\rangle\leq\mu_i(\omega)$, since $\omega$ is strictly dominant. 
\end{remark}

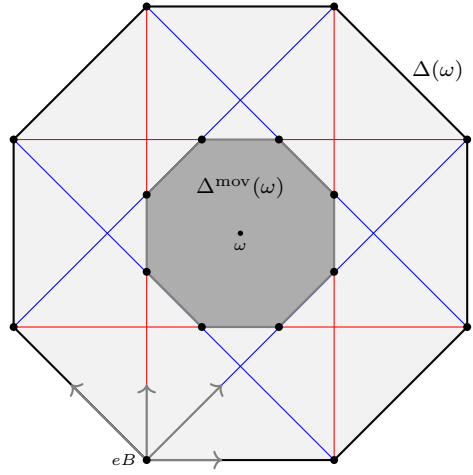
\begin{figure}[ht]
    \centering    
\begin{tikzpicture}
    \draw[thick] (1.24,3) -- (3,1.24) -- (3,-1.24) -- (1.24,-3) -- (-1.24,-3) -- (-3,-1.24) -- (-3,1.24)-- (-1.24,3) -- cycle;
    \fill[gray, opacity = 0.1] (1.24,3) -- (3,1.24) -- (3,-1.24) -- (1.24,-3) -- (-1.24,-3) -- (-3,-1.24) -- (-3,1.24)-- (-1.24,3) -- cycle;
    \node[right, font = \footnotesize] at (2.15,2.15) {$\Delta(\omega)$};
 
    \draw[red] (-1.24,-3) -- (-1.24,3);
    \draw[red] (1.24,-3) -- (1.24,3);
    \draw[red] (-3,-1.24) -- (3,-1.24);
    \draw[red] (-3,1.24) -- (3,1.24);
    \draw[blue] (-1.24,-3) -- (3,1.24);
    \draw[blue] (1.24,-3) -- (-3,1.24);
    \draw[blue] (1.24,3) -- (-3,-1.24);
    \draw[blue] (-1.24,3) -- (3,-1.24);
    
   \draw[thick, gray] (1.24,0.51) -- (0.51,1.24) -- (-0.51,1.24) -- (-1.24,0.51) -- (-1.24,-0.51) -- (-0.51,-1.24) -- (0.51,-1.24) -- (1.24,-0.51) -- cycle;
   \fill[gray, opacity = 0.6] (1.24,0.51) -- (0.51,1.24) -- (-0.51,1.24) -- (-1.24,0.51) -- (-1.24,-0.51) -- (-0.51,-1.24) -- (0.51,-1.24) -- (1.24,-0.51) -- cycle;
   \node[font = \footnotesize] at (0,0.6) {$\Delta^{\mov}(\omega)$};
   
   \fill (-0.51,1.24) circle[radius = 1.5pt];
   \fill (0.51,1.24) circle[radius = 1.5pt];
   \fill (0.51,-1.24) circle[radius = 1.5pt];
   \fill (-0.51,-1.24) circle[radius = 1.5pt];    

   \fill (1.24,0.51) circle[radius = 1.5pt];
   \fill (1.24,-0.51) circle[radius = 1.5pt];
   \fill (-1.24,-0.51) circle[radius = 1.5pt];  
   \fill (-1.24,0.51) circle[radius = 1.5pt];
   
    \fill (-3,1.24) circle[radius = 1.5pt];
   \fill (3,1.24) circle[radius = 1.5pt];
   \fill (-3,-1.24) circle[radius = 1.5pt];
   \fill (3,-1.24) circle[radius = 1.5pt];
   \fill (1.24,3) circle[radius = 1.5pt];
   \fill (1.24,-3) circle[radius = 1.5pt];
   \fill (-1.24,3) circle[radius = 1.5pt];

   \fill (0,0) circle[radius = 1pt];
   \node[below, font = \tiny] at (0,0) {$\omega$};
   
    \draw[->,  thick, gray](-1.24,-3) -- (-0.24,-3);
    \draw[->,  thick, gray](-1.24,-3) -- (-2.24,-2);
    \draw[->,  thick, gray](-1.24,-3) -- (-1.24,-2);
    \draw[->,  thick, gray](-1.24,-3) -- (-0.24,-2);
    \fill (-1.24,-3) circle[radius = 1.5pt];   
    \node[ left, font = \tiny] at (-1.25,-3) {$eB$};    
\end{tikzpicture}
    \caption{\centering The polytopes $\Delta(\omega)$ and $\Delta^{\mov}(\omega)$ for $\omega$ in $\DC_2$.}
     \label{fig:polytopes}
\end{figure}

For the sake of clarity, we illustrate in Figure \ref{fig:polytopes} both the polytope of weights and of movable weights for a choice of $\omega$ in the type  $\DC_2$ case, with their GIT chamber decomposition. The origin, corresponding to $eB$, is decorated with the set of positive roots $\Phi^+$.

\section{Existence of MG quotients}\label{sec:mainproof}

In this section we address the proof of the main result of the paper, which we have introduced in Theorem \ref{theorem:Introduction} but we state below with full details.
\begin{theorem}\label{thm:main}
In the situation of Setup \ref{set:flags}, there exists a geometric quotient $\MX$ of $X$ such that:
\begin{itemize}[leftmargin=\marit pt]
\item $\MX$ is an MDS, and there is an isomorphism of $\C$-algebras $\Cox(\MX)\simeq\Cox(X)$;
\item $\MX$ is an MG quotient of $X$;
\item the birational map $\Env(\MX)\to \CX$ is an isomorphism. 
\end{itemize}
\end{theorem}

\begin{remark}\label{rem:main}
Let us briefly discuss the two cases we have excluded in the Setup \ref{set:flags}. In the case of $\DA_1$, the associated flag variety is $\P^1$; all its GIT quotients are isomorphic to a point, hence all of them are MG quotients, and we have $\Env(\MX)\simeq \CX$. In the case of $\DA_2$, every geometric quotient of $X$ is isomorphic to $\P^1$, hence every geometric quotient is an MG quotient, and $\Env(\MX)\simeq\CX\simeq \P^1$. On the other hand, it is not true in these cases that $\Cox(X)\simeq \Cox(\MX)$.
\end{remark}

Throughout this section we will work in the situation of Setup \ref{set:flags}. We will show that, given a general element $v\in \Delta^{\mov}({{\omega}})$ for a good weight $\omega$, the GIT quotient $\GX_v$ is an MG quotient satisfying the requirements of Theorem \ref{thm:main}. Here the word {\em general} in a polytope (resp. cone) means ``in the complement of a finite union of polytopes (resp. cones) of smaller dimension''. We start by noting the following known relevant facts.

\begin{remark}\label{rem:Qfact}
In the above situation, if $v\in \Delta^{\mov}({{\omega}})\cap \Mo(H)_\Q$ is general, then we may assume that $\GX_v$ is a geometric quotient of $X$. Then,  by \cite[Lemma~2.1]{HuKeel}, $\GX_v$ is $\Q$-factorial and, subsequently, an MDS by \cite[Theorem~1.1]{Baker}.
\end{remark}

The key idea for our proof is that the general point of every fundamental divisor $D_i^\sigma$ is semistable with respect to the chosen linearization $v$. We will denote by $X^\ss_v\subset X$ the set of $v$-semistable points for a sufficiently large power of $L({{\omega}})$.

\begin{proposition}\label{prop:Pmovseesalldiv}
In the situation of Setup \ref{set:flags}, let $\omega\in\Q_{>0}(\omega_1,\dots,\omega_n)$ be a good weight, and $v\in \Delta(\omega)\cap \Mo(H)_\Q$ be a general element. Then $v\in \Delta^{\mov}(\omega)$ if and only if $\codim (X\setminus X^\ss_v)\geq 2$. In this case, for every $i,\sigma$, the open subset $X^\ss_v\cap D_i^\sigma$ is dense in $D_i^\sigma$ and its quotient is a prime divisor in $\GX_v$. 
\end{proposition}

\begin{proof}
Since $v$ is general, we assume that $\GX_v$ is a geometric quotient of $X$. 
By Hilbert--Mumford criterion, the unstable locus in $X$ for the linearization $v$ is the union of the unstable loci of every $\mu(v)$, $\mu\in\No(H)$. By \cite{BBS} (see also \cite[Section~2.2]{WORS3}), the complement of the (semi)stable locus $X^\ss_{\mu(v)}$ of the $\C^*$-action induced by $\mu$ for the linearization $\mu(v)\in\Q$ is a union of (closures of) BB-cells. 

Now, by using Lemma \ref{lem:omegaifundi}, Remark \ref{rem:omegaifundi}, and Corollary \ref{cor:BBcellsfund}, the only possible divisors that may appear in the unstable locus of $v$ are the fundamental divisors. As in Lemma \ref{lem:omegaifundi}, we may write every $D_i^\sigma=\ol{X^-(Y_i^\sigma)}$, where $Y^\sigma_i$ is a $\mu_i^\sigma$-fixed point component, whose $L(\omega)$-weight is, by Remark \ref{rem:omegaifundipt}, equal to $\mu_i^\sigma(\omega)-\mu_i^{r_i}(\omega)$.

Summing up, $X\setminus X^\ss_{v}$ contains a divisor if and only if there exists $i,\sigma$ such that $\mu_i^\sigma(v)\leq \mu_i^\sigma\mu_{L(\omega)}(Y_i^\sigma)=\mu_i^\sigma(\omega)-\mu_i^{r_i}(\omega)$, i.e., if and only if $v\notin \Delta^{\mov}(\omega)$.
\end{proof}

\begin{lemma}\label{lem:isotropy}
In the situation of Proposition \ref{prop:Pmovseesalldiv}, $H$ acts freely on an $H$-invariant open subset $O\subseteq X$ whose complement has codimension at least two.  
\end{lemma}

\begin{proof}
For every $\sigma \in W$ we will consider the Bruhat open neighborhood $U_\sigma\subset X$ of $\sigma B\in X$, introduced in Section \ref{ssec:CoxRH}. 

Let $O\subset X$ be the $H$-invariant open subset of points with trivial isotropy. Since the $U_\sigma$'s form an open covering of $X$, it is enough to show that the closed sets $(X\setminus O)\cap U_\sigma$ have codimension at least two in $U_\sigma$ for every $\sigma$ and, by considering the action of $W$, it is enough to show this for $\sigma=e$. 

Identifying $U_e\simeq\fn$, the set $(X\setminus O)\cap U_e$ corresponds to a union $Z$ of coordinate vector subspaces, namely
\[
Z=\{x\in \fn\mid \Z Q(x)\subsetneq \Mo(H)\},
\]
where $Q(x):=\{\alpha\in \Phi^+\mid x_{\alpha}\neq 0\}$,  for $x =\textstyle \sum_{\alpha\in \Phi^+}x_\alpha$ with each  $x_\alpha\in \fg_{-\alpha}$. Then the statement reduces to show that, for every $\alpha\in\Phi^+$, the set $\Phi^+\setminus\{\alpha\}$ generates $\Mo(H)$, which is obvious since the Dynkin diagram of $\fg$ is connected.
\end{proof}

\begin{proposition}\label{prop:Picquotient}
In the situation of Proposition \ref{prop:Pmovseesalldiv}, the pullback map induces an isomorphism $\Cl(\GX_v)\simeq \Pic^H(X)$.
\end{proposition}

\begin{proof}
With the notation of Lemma \ref{lem:isotropy}, set $O_v:=O\cap X^{\ss}_v$, whose complement has codimension at least two by Proposition \ref{prop:Pmovseesalldiv} and the Lemma \ref{lem:isotropy}. The geometric quotient $\cG O_v$ of the $H$-invariant open subset $O_v$ is a smooth open subset of $\GX_v$, whose complement has still codimension at least two. Hence the quotient $O_v\to \cG O_v$ gives the following pullback map
\begin{equation}\label{eq:pullback}
\Cl(\GX_v)\simeq\Cl(\cG O_v)\simeq\Pic(\cG O_v)\longrightarrow \Pic^H(O_v)\simeq\Pic^H(X),
\end{equation}
which is known to be injective. 
Then, given a linearized line bundle $L\in \Pic^H(X)$, \cite[Th\'eor\`eme~2.3]{DN89}  tells us that its restriction to $O_v$ is a pullback of a line bundle in $\cG O_v$ if and only if the stabilizer $H_x$ acts trivially on the fiber $L_x$, for every $x\in O_v$. This is obvious by Lemma \ref{lem:isotropy}, thus we conclude that (\ref{eq:pullback}) is an isomorphism. Note that the sheaf in $\GX_v$ corresponding to the descent of a linearized line bundle $L$ is, in general, only a rank one reflexive sheaf.
\end{proof}

\begin{proposition}\label{prop:Coxquotient}
In the situation of Proposition \ref{prop:Pmovseesalldiv}, $\Cox(\GX_v)\simeq \Cox(X)$ as $\C$-algebras. Moreover, $\GX_v$ is an MG quotient of $X$ and $\Env(\GX_v)\simeq \CX$.
\end{proposition}

\begin{proof}
With the notation of the above proposition, let $L(\lambda,m)\in\Pic(\cG O_v)$ be the line bundle corresponding to the pair $(\lambda,m)\in\Pic^H(X)$ via the pullback isomorphism induced by the quotient map. In other terms, $L(\lambda,m)$ is the pushforward of the subbundle of $H$-invariant sections of $L(\lambda)|_{O_v}$ with respect to the linearization $m$. Therefore, there is a natural isomorphism $\HH^0(X, L(\lambda))_m\simeq\HH^0(\cG X_v, L(\lambda,m))$ applying Hartogs' extension theorem. These isomorphisms are compatible with the multiplication of sections, and hence provide an isomorphism of $\Pic^H(X)$-graded algebras, where $\Cox(X)$ is endowed with the finer grading introduced in Section~\ref{sec:MGRH}. In fact, choosing $\{(\omega_i,0),(0,\alpha_i)\mid i\in\{1,\dots,n\}\}$ as a $\Z$-basis of $\Pic^H(X)\simeq\Cl(\GX_v)$, we have
\begin{align*}
\Cox(\GX_v)&=\bigoplus_{a_1,\dots,a_n\in\Z}\bigoplus_{b_1,\dots,b_n\in\Z}{\textstyle\HH^0\big(\GX_v,L\big(\sum_{i=1}^n a_i\omega_i,\sum_{i=1}^n b_i\alpha_i\big)\big)}\\&\simeq\bigoplus_{a_1,\dots,a_n\in\Z}\bigoplus_{b_1,\dots,b_n\in\Z}{\textstyle\HH^0\big(X,L\big(\sum_{i=1}^n a_i\omega_i\big)\big)_{\sum_{i=1}^n b_i\alpha_i}}\\&=\bigoplus_{a_1,\dots,a_n\in\Z}{\textstyle\HH^0\big(X,L\big(\sum_{i=1}^n a_i\omega_i\big)\big)}=\Cox(X).
\end{align*}

For the second part it is enough to note that from the above proof it follows that every GIT quotient of $X$, which is of the form $\Proj (\bigoplus_{k\geq 0}\HH^0(X,L(k\lambda))_{km})$ for some $(\lambda,m)$, can be written as $\Proj( \bigoplus_{k\geq 0}\HH^0(\GX_v,L(k\lambda,km)))$, and so it is a rational contraction of $\GX_v$; hence we conclude that $\GX_v$ is an MG quotient. Conversely, the same argument tells us also that every rational contraction of $\GX_v$ is a GIT quotient of $X$, hence $\Env(\GX_v)\simeq \CX$.
\end{proof}

At this point, to conclude the proof of Theorem \ref{thm:main}, we just need to verify that every $X$ has at least a good weight.
\begin{proposition}\label{prop:Deltamovnotempty}
The weight $\ol\omega:=\sum_{i=1}^n\omega_i$ is good.
\end{proposition}

\begin{proof}
It is enough to show that, for every $i,\sigma$, the open half-space defined by the inequality $\mu_i^\sigma(m)> \mu_i^\sigma({\ol{\omega}})-\mu_i({\ol{\omega}})+1$ contains ${\ol{\omega}}$. Then it suffices to show that the components of ${\ol{\omega}}$ with respect to the simple roots are all larger than one. This can be readily checked, for instance, thanks to the formulas presented in \cite[p.~69]{Hum}.  
\end{proof}

\begin{remark}\label{rem:Deltamovgen} 
From the proof of Proposition \ref{prop:Pmovseesalldiv}, we also see that the geometric quotients of $X$ coming from linearizations in $\Delta({\omega})\setminus \Delta^{\mov}(\omega)$ provide rational contractions of any MG quotient $\MX$ of $X$, in which at least one of the divisors $D_i^\sigma$ gets contracted. Summing up, we may conclude the following identities:
\begin{align*}\label{eq:effmov}
\Eff(\MX)&=\left\{(\omega,m)\in \Pic^H(X)_\R\mid \mu_i^\sigma(m)\geq \mu_i^\sigma(\omega)-\mu_i(\omega),\mbox{ for every }i,\sigma\right\}=\Delta,\\
\Mov(\MX)&=\left\{(\omega,m)\in \Pic^H(X)_\R\mid \mu_i^\sigma(m)\geq \mu_i^\sigma(\omega)-\mu_i^{r_i}(\omega),\mbox{ for every }i,\sigma\right\}\subset\Delta.
\end{align*}
\end{remark}

\section{Construction of an MG quotient}\label{sec:ConsMG}

In this section we explain how to construct an MG quotient of the complete flag variety $X = G/B$. To this end, we first describe MG quotients of local affine patches and study their contractions of relative Picard number $1$. We then discuss how to obtain an MG quotient $\MX$ from these data.

As in the previous sections, we assume $G$ to be an adjoint simple algebraic group and $B$ a Borel subgroup containing the maximal torus $H$. Denote by $\Pi:=\{\alpha_1,\dots,\alpha_n\}\subset\Phi^+$ the corresponding positive simple roots, with $n:=\rk(H)\geq 2$ and $r:=|\Phi^+-\Pi|$. We always assume that the Dynkin diagram of $G$ is neither $\DA_1$ nor $\DA_2$.

\begin{notation}\label{not:Disigma}
For every $i=1,\ldots,n$ and $\sigma\in W$, we will abuse notation by denoting by \[D^\sigma_i\subset\GX_v\]  the quotient of the fundamental divisor $D^\sigma_i$ of $X$ in any GIT quotient $\GX_v$ of $(X,L(\omega))$, with $v$ satisfying $\mu_i^{\sigma}(v)\geq \mu_i^{\sigma}(\omega)-\mu_i^{r_i}(\omega)$; note that these are precisely the GIT quotients whose sets of semistable points $X^{\ss}_v\subset X$ intersect the fundamental divisor $D^\sigma_i\subset X$. In particular, $D_i^\sigma$ represents a divisor in every MG quotient of $X$.
\end{notation}

\begin{remark}\label{rem:MXuponMU}
The idea for constructing MG quotients of $X$ that we propose in this section essentially comes from the following observations.
\begin{itemize}[leftmargin=\marit pt]
\item By \cite[Proposition~3.8]{OS26}, we can realize the geometric quotients of $U_e$ as geometric quotients of $X$. More concretely, given an ample line bundle $L(\omega)$ in $X$, the geometric quotients of $U_e$ correspond to the GIT chambers of $\Delta(\omega)$ that intersect every neighborhood of the weight $0\in \Mo(H)_\R$. Considering the convex cones generated by these chambers, we obtain the GIT chamber decomposition of the cone $\R_{\geq0}\Phi^+$, which is precisely the GIT chamber decomposition for the $H$-action on $U_e\simeq\fn$.

The next section will describe how to effectively construct these MG quotients of $U_e$.
\item The complement $X\setminus U_e$ is the union of all the fundamental divisors $D_i^e$. This means that an MG quotient of $U_e$ must come from a stable locus that meets every $D_i^\sigma$, $\sigma\neq e$; moreover, in order to construct an MG quotient $\MX$ of $X$ upon an MG quotient of $U_e$, we must perform $n$ divisorial extractions centered in the sets parametrizing $H$-orbits contained in the BB-cells $\ol{X^+(Y_i^e)}$ (with the notation of Lemma \ref{lem:omegaifundi}), producing exactly the divisors $D^e_i\subset\MX$. In principle, we may need to perform these extractions in a particular order, and to perform also some small $\Q$-factorial modifications. However, we will show that the MG quotient of $U_e$ can be chosen so that we do not need any extra small $\Q$-factorial modifications and the divisorial extractions can be performed in any order. 
\end{itemize}
\end{remark}

\begin{figure}[ht!]
	\centering    
		\begin{tikzpicture}
			\draw[thick] (1.24,3) -- (3,1.24) -- (3,-1.24) -- (1.24,-3) -- (-1.24,-3) -- (-3,-1.24) -- (-3,1.24)-- (-1.24,3) -- cycle;
			\fill[gray, opacity = 0.1] (1.24,3) -- (3,1.24) -- (3,-1.24) -- (1.24,-3) -- (-1.24,-3) -- (-3,-1.24) -- (-3,1.24)-- (-1.24,3) -- cycle;

			\fill[gray, opacity = 0.3] (-1.24,-3) -- (0,-1.75) -- (-0.51,-1.24) -- (-1.24,-1.24) -- cycle;
			\node[font = \tiny] at (-0.82,-2) {$\MU_e$};
			
			\fill[gray, opacity = 0.6] (-1.24,-1.24) -- (0.51,-1.24) -- (0,-1.75) -- (-1.24,-0.51) -- cycle;
			
			\fill[gray, opacity = 0.8] (1.24,0.51) -- (0.51,1.24) -- (-0.51,1.24) -- (-1.24,0.51) -- (-1.24,-0.51) -- (-0.51,-1.24) -- (0.51,-1.24) -- (1.24,-0.51) -- cycle;

			\draw[red] (-1.24,-3) -- (-1.24,3);
			\draw[red] (1.24,-3) -- (1.24,3);

			\draw[red] (-3,-1.24) -- (3,-1.24);
			\node[above, red,  font = \tiny] at (-2, -1.3) {$D^e_2$};
			\draw[red] (-3,1.24) -- (3,1.24);

			\draw[blue] (-1.24,-3) -- (3,1.24);
			\draw[blue] (1.24,-3) -- (-3,1.24);
			\node[above, rotate = -45, blue, font = \tiny] at (0.5,-2.3) {$D^e_1$};
			\draw[blue] (1.24,3) -- (-3,-1.24);
			\draw[blue] (-1.24,3) -- (3,-1.24);
			
			\fill (-3,1.24) circle[radius = 1.5pt];
			\fill (3,1.24) circle[radius = 1.5pt];
			\fill (-3,-1.24) circle[radius = 1.5pt];
			\node[ left, font = \tiny, black] at (-3,-1.24) {$r_2B$};
			\fill (3,-1.24) circle[radius = 1.5pt];
			\fill (1.24,3) circle[radius = 1.5pt];
			\fill (1.24,-3) circle[radius = 1.5pt];
			\node[ right, font = \tiny, black] at (1.25,-3) {$r_1B$};
			\fill (-1.24,3) circle[radius = 1.5pt];
			\fill (-1.24,-3) circle[radius = 1.5pt];   
			\node[ left, font = \tiny, black] at (-1.25,-3) {$eB$};
			
			\fill (-0.51,1.24) circle[radius = 1.5pt];
			\fill (0.51,1.24) circle[radius = 1.5pt];
			\fill (0.51,-1.24) circle[radius = 1.5pt];
			\fill (-0.51,-1.24) circle[radius = 1.5pt];    
			
			\fill (1.24,0.51) circle[radius = 1.5pt];
			\fill (1.24,-0.51) circle[radius = 1.5pt];
			\fill (-1.24,-0.51) circle[radius = 1.5pt];  
			\fill (-1.24,0.51) circle[radius = 1.5pt];
			
			\draw[->, red, opacity = 0.7] (-0.72,-1.74) to [bend left = 47] (-1,-1);
			\draw[->, blue, opacity = 0.7] (-1,-1) to [bend left = 43] (-0.31,-0.74);
			\draw[->, blue, opacity = 0.7] (-0.72,-1.74) to [bend right = 47] (0,-1.41);
			\draw[->, red, opacity = 0.7] (0,-1.41) to [bend right = 43] (-0.31,-0.74);
			
			\node[font = \tiny] at (0,0) {$\MX$};			
		\end{tikzpicture}
	\caption{\centering How to reach $\MX$ from $\MU_e$.}
	\label{fig:sec5}
\end{figure}
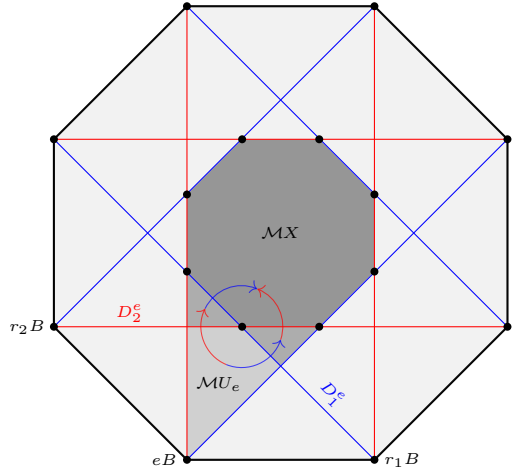

We may intuitively represent the second step of the above Remark by means of Figure \ref{fig:sec5} for $n=2$; starting from an appropriate MG quotient of $U_e$, we can obtain an MG quotient of $X$ extracting, in any order, the $n$ fundamental divisors $D^e_1,\ldots,D_n^e$. In terms of GIT chamber decomposition, this corresponds to the existence of paths connecting a chamber of $\Delta(\omega)$ containing the origin to a chamber in $\Delta^{\mov}(\omega)$, which cross exactly $n$ divisorial walls.

\subsection{ MG quotients of Bruhat neighborhoods and their birational geometry}\label{ssec:MGUe} 

We consider the $H$-invariant Bruhat open neighborhood $U_e$  of $eB$ in $X$, $H$-equivariantly isomorphic to  $\fn=\bigoplus_{\alpha\in \Phi^+}\fg_{-\alpha}$ (see Section \ref{ssec:CoxRH}). 
Note that the Bruhat decomposition, translated into our language, provides
\[
X=U_e\sqcup \bigcup_{i=1}^nD_i^e,
\]
so that GIT quotients of $U_e$ can be seen as GIT quotients $\GX_v$ of polarized pairs $(X,L(\omega))$, with $L(\omega)$ ample, such that $v\in \Mo(H)_\Q$ satisfies that
\[v\in \Delta(\omega)\cap\{m\mid \mu_i(m)\leq \langle \omega,\alpha_i\rangle,\mbox{ for every } i\}.
\]
In fact, $v$ satisfies the required extra inequalities if and only if $X^\ss_v\cap D_i^e=\emptyset$ for every $i$, that is, if and only if $X^\ss_v\subseteq U_e$. This shows that every quotient $\GX_v$ is also a GIT quotient of $U_e$. The converse follows from the proof of \cite[Proposition~3.8]{OS26}. 

On the other hand $U_e$ is a toric variety isomorphic, via the exponential map, to the complex affine space $\fn\simeq\C^{n+r}$, whose cone is the dual of the non-negative orthant $\sigma^+\subset\Mo(T)_\R\simeq\R^{n+r}$, where $T:= \prod_{\alpha\in \Phi^+}(\fg_{-\alpha}\setminus\{0\})$. 
Then the GIT quotients of $U_e$ with respect to the induced $H$-action are in correspondence with the chambers of the GIT chamber decomposition of $\pi_\R(\sigma^+)$, where $\pi : \Mo(T)\to\Mo(H)$ is the natural {\it weight map}. In suitable coordinates, $\pi$ is given by a {\it weight matrix} 
\begin{equation*}\label{eq:weightmatrix}
	A:=\left( \,{I} _n \mid  M\,\right), \quad  M\in \Z_{\geq0}^{n\times r},
\end{equation*}
where the columns of $ M$ are the coordinates of the $r$ positive non-simple roots with respect to the fixed basis of simple roots, in some order.

Given any $v\in\pi_\Q(\sigma^+):=\pi_\R(\sigma^+)\cap\Mo(H)_\Q$, the GIT quotient $\cG U_{e,v}$ is the $(T/H)$-toric variety whose moment polytope is $P_v := \pi_\R^{-1}(v)\cap\sigma^+$, which is bounded by \cite[Lemma 3.2]{OS26}. The polytope $P_v$ has then a nice explicit description, which only depends on the positive non-simple roots and the choice of $v$. 

\begin{lemma}
Under the natural identification $\ker(\pi_\R)\simeq\Mo(T/H)_\R\simeq\R^r$, the moment polytope $P_v$ of the GIT quotient $\cG U_{e,v}$ is 
\begin{equation}\label{eq:P_v}
	P_v = \left\{x \in\R^r_{\geq0} \mid - M x+ v\geq  0\right\},
\end{equation}
where $v$ denotes the column matrix whose coordinates are $v_i=\mu_i(v)$, $i=1,\dots,n$.
\end{lemma}

\begin{proof}
	The kernel $\ker(\pi_\R)$ is generated in $\Mo(T)_\R$ by the columns of the matrix 
	\begin{equation*}\label{eq:galedual}
	A^\perp:=\begin{pmatrix}-M\\\hline I_r\end{pmatrix}. 
\end{equation*}
Considering the element $v'$ in $\pi_\R^{-1}(v)\subseteq\Mo(T)_\R\simeq \R^{n+r}$ whose coordinates are $(v_1,\dots,v_n,0,\dots,0)$, we then have that $P_v=\{v'+A^\perp x \mid x\in \R^r\}\cap\sigma^+$, from which
the identity (\ref{eq:P_v}) follows.
\end{proof}

Throughout the rest of the section we will consider the {\it Gale transform} \[\pi^\perp_\R:\No(T)_\R\to \No(T/H)_\R\] of $\pi_\R:\Mo(T)_\R\to \Mo(H)_\R$, which is represented by the transpose $\left(- M^t\mid  I_r\right)$ of the matrix $A^\perp$ introduced in the above proof. 
We will denote the first $n$ columns of this matrix by \[-\nu_1,\dots,-\nu_n.\] The last $r$ columns are in one-to-one correspondence with positive non-simple roots (in the order we have chosen to define $A$ above), so we will denote them by \[e_\alpha, \,\, \alpha \in \Phi^+\setminus \Pi.\]

Our next goal will be to find an MG quotient of $U_e$. Whilst the existence of geometric quotients of an affine space satisfying the required properties can be inferred from \cite[Proposition 15.2.4]{CLS}, we will explicitly write an argument that allows us to describe precisely which geometric quotients of $U_e$ are MG quotients. In the next statements we will use the following notation:
 \[K:=\R_{\geq0}(\Phi^+\setminus\Pi)\subset\Mo(H)_\R.\]

\begin{proposition}\label{prop:MUe}
    Let $v\in\pi_\Q(\sigma^+)$ be such that $\GU_{e,v}$ is a geometric quotient. 
    Then $\GU_{e,v}$ is an MG quotient of $U_e$ if and only if $v\in K$. In this case $\MU_e:=\GU_{e,v}$ is an $r$-dimensional toric variety whose fan rays are generated by $-\nu_1,\dots,-\nu_n$ and every  $e_\alpha,\,\,\alpha \in \Phi^+\setminus \Pi$. 
    Moreover, there is an isomorphism $\NU(\MU_e)\simeq\Mo(H)_\R$, that identifies the movable cone of $\MU_e$ with the cone generated by the positive non-simple roots of $\fg$, that is 
\begin{equation}\label{eq:mov}
	\Mov(\MU_e)= \R_{\geq0}(\Phi^+\setminus\Pi).
\end{equation}	
\end{proposition}

\begin{proof}
Note that, in the situation of Setup \ref{set:flags}, the cone $K$ has dimension $n$. 

Let us consider a general element $v\in \Int(K)\cap \Mo(H)_\Q$; it can be written as a positive rational combination of all the positive non-simple roots, in other words, there exists $x_0\in \Q_{>0}^r$ such that 
\[Mx_0=v\in \Int(K).\] 
We will now describe the facets of the polytope $P_v$, whose normals define the rays of the associated normal fan. The fact that it has $r$ facets normal to $e_\alpha, \,\, \alpha \in \Phi^+\setminus \Pi$ follows since $P_v$  contains elements of the form $\varepsilon e_\alpha$, $0<\varepsilon\ll 1$ for every $\alpha\in \Phi^+\setminus \Pi$. On the other hand, since $K$ is $n$-dimensional, for every $i$ there exists $y(i)\in \Q_{>0}^r$ such that $v(i):=My(i)$ has coordinates $\mu_i(v(i))=\mu_i(v)$, and $\mu_j(v(i))<\mu_j(v)$ for $j\neq i$. We conclude that $P_v$ has $n$ facets normal to $-\nu_1,\dots,-\nu_n$.

On the other hand, the arguments above tell us that the set of facet normals of a polytope $P_v'$ associated to any other $v'\in \pi_\Q(\sigma^+)$ is contained in $\{e_\alpha\mid  \alpha \in \Phi^+\setminus \Pi\}\cup\{-\nu_1,\dots,-\nu_n\}$, from which it follows that any other $\GX_{v'}$ is a rational contraction of our $\GX_v$. 
Note that if $v$ is general (so that it lies outside the inner walls of the GIT chamber decomposition), $\GX_v$ is also $\Q$-factorial. Summing up, we conclude that $\MU_e:=\GX_v$ is an MG quotient of $U_e$. 

The fact that $\NU(\MU_e)\simeq\Mo(H)_\R$ follows then from our description of the rays of the fan of $\MU_e$ and standard toric geometry. It is also a standard fact that under this identification $\Eff(\MU_e)$ coincides with $\pi_\R(\sigma^+)=\R_{\geq 0}\{\alpha\mid \alpha\in \Phi^+\}$. We will now prove that, via that isomorphism, $K$ corresponds to the movable cone of $\MU_e$. Note that our previous argument tells us that all the $\GX_v$, for $v$ general in $K$, are isomorphic in codimension one, hence $K\subseteq \Mov(\MU_e)$. 

In order to prove that $\Mov(\MU_e)\subseteq K$, let us first note that $K$ can be also described as \[K= \big\{m\in\Mo(H)_\R\mid \mu_i(m)\geq0,\,\mu_i^{r_i}(m)\geq0,\mbox{ for every }i\big\}.\]
A proof of this fact can be found in \cite[Theorem~4.2]{Pany}. 

We now write
\[\pi_\R(\sigma^+)=K\cup\bigcup_{i=1}^n N^i, \textrm{ where } N^i:=\big\{m\mid \mu_j(m)\geq0 \mbox{ for every }j,
\,\,\mu_i^{r_i}(m)<0\big\}.\]
The proof may be concluded by showing that, for every $i$, an element $u\in N^i$ satisfies that the facet normals of $P_u$ are contained in $\{e_\alpha\mid\alpha \in \Phi^+\setminus \Pi\}\cup\{-\nu_j\mid j\neq i\}$. Note that this also shows that weights in $\pi_\R(\sigma^+)\setminus K$ do not provide MG quotients of the variety $U_e$.

Given $u\in N^i$ we have that $\mu_i^{r_i}(u)<0$, 
which can be written as:
\begin{equation}\label{eq:inequ}
\mu_i(u)>-\sum_{j\neq i}\langle \alpha_j,\alpha_i\rangle\mu_j(u).\end{equation} Let $x$ be an element satisfying 
\begin{equation}\label{eq:Punoti}
x\geq 0,\quad\mbox{and}\quad \nu_j x\leq \mu_j(u),\mbox{ for every }j\neq i;
\end{equation}
We will show that then $\nu_i x\leq \mu_i(u)$, i.e., that $x\in P_u$. In other words, this shows that $-\nu_i$ is not a facet normal of $P_u$. 

Now if $x$ satisfies the inequalities (\ref{eq:Punoti}), set $y=Mx$, so that we have $\mu_j(y)\leq \mu_j(u)$, for $j\neq i$, and we want to show that the inequality holds also for $j=i$. 

Since $x\geq 0$, then $y\in K$, and so 
\begin{equation}\label{eq:muiri}
0\leq \mu_i^{r_i}(y)=-\mu_i(y)-\sum_{j\neq i}\langle \alpha_j,\alpha_i\rangle\mu_j(y),
\end{equation}
so we have $\mu_i(y)\leq -\sum_{j\neq i}\langle \alpha_j,\alpha_i\rangle\mu_j(y)$. Since $\langle \alpha_j,\alpha_i\rangle\leq 0$ and $\mu_j(y)\leq \mu_j(u)$ for $j\neq i$, and taking into account the inequality (\ref{eq:inequ}), we finally get
\[\mu_i(y)\leq -\sum_{j\neq i}\langle \alpha_j,\alpha_i\rangle\mu_j(y)\leq -\sum_{j\neq i}\langle \alpha_j,\alpha_i\rangle\mu_j(u)<\mu_i(u).\qedhere\]
\end{proof}

\begin{remark}\label{rem:Cl(MUe)}
	The divisor class group of $\MU_e$ is identified with $\Mo(H)$ via the weight map $\pi$ itself. This implies that, for every $i\in\{1,\dots,n\}$, the $(T/H)$-invariant divisor corresponding to the negative ray $-\nu_i$ has numerical class $\alpha_i\in \Mo(H)_\R\simeq \NU(\MU_e)$. On the other hand, one may easily see that this divisor is the image of $D_i^{r_i}\subset X$, which we simply denote $D_i^{r_i}\subset \MU_e$. Moreover, the next proposition will tell us that crossing the wall separating $K$ and $N^i$, defined by $\mu_i^{r_i}(m)=0$, corresponds to a rational divisorial contraction of $D_i^{r_i}\subset \MU_e$. 
\end{remark}

\begin{proposition}\label{prop:foundcontraction}
With the same notation as above, for every $i\in\{1,\dots,n\}$ set
\[ K^i:=\big\{m\mid \mu_j(m)\geq0 \mbox{ for every }j,\,\,\mu_j^{r_j}(m)\geq0,\mbox{ for every }j\neq i, 
\,\,\mu_i^{r_i}(m)<0\big\}.\]
Then, for every $u\in\Int(K^i)$ the set of facet normals of $P_u$ is equal to 
\[\{e_\alpha\mid\alpha \in \Phi^+\setminus \Pi\}\cup\{-\nu_j\mid j\neq i\}.\] 
In particular, given an MG quotient $\MU_e:=\GU_{e,v}$, $v\in \Int(K)\cap \Mo(H)_\Q $, and a general element $u\in\Int(K^i)\cap \Mo(H)_\Q$, the rational map $\MU_e\dashrightarrow \GU_{e,u}$ is the composition of an SQM and the divisorial contraction of the strict transform of $D_i^{r_i}\subset \MU_e$. 
\end{proposition}

\begin{figure}[ht]
	\centering    
		\begin{tikzpicture}[scale=1.2, x={(1cm,0cm)}, y={(0.5cm,0.8660254cm)}]
			
			\draw[thick] (0,0) -- (0,6) -- (6,0) -- cycle;
			\fill[gray, opacity = 0.3] (0,3) -- (3,3) -- (4,2) -- (3,1.5) -- cycle;
			\fill[red, opacity = 0.25] (0,0) -- (0,3) -- (3,1.5)  -- cycle;
			\fill[green, opacity = 0.25] (0,6) -- (0,3) -- (3,3)  -- cycle;
			\fill[blue, opacity = 0.25] (6,0) -- (3,1.5) -- (4,2)  -- cycle;
			\draw[very thin] (0,0) -- (4,2);
			\draw[very thin] (0,3) -- (6,0);
			\draw[very thin] (0,3) -- (3,3);
			\draw[very thin, dashed, opacity=0.5] (3,1.5) --(0,6);
			\draw[very thin, dashed, opacity=0.5] (2,2) --(0,6);
			\draw[very thin, dashed, opacity=0.5] (3,1.5) --(3,3);
			\draw[very thin, dashed, opacity=0.5] (2,2) --(4,2);
			\draw[very thin, dashed, opacity=0.5] (0,3) --(4,2);
			\draw[very thin, dashed, opacity=0.5] (0,0) -- (3,3);
			\draw[very thin, dashed, opacity=0.5] (2.4,2.4) --(6,0);
			\node[ left, font = \tiny, black] at (0,0) {$\alpha_1$};
			\node[ left, font = \tiny, black] at (0,6) {$\alpha_2$};
			\node[ right, font = \tiny, black] at (6,0) {$\alpha_3$};
			\node[ right, font = \tiny, black] at (3,3) {$\alpha_2+\alpha_3$};
			\node[ right, font = \tiny, black] at (4,2) {$\alpha_2+2\alpha_3$};
			\node[ left, font = \tiny, black] at (0,3) {$\alpha_1+\alpha_2$};
			\node[ left, font = \tiny, black] at (2,2) {$\alpha_1+\alpha_2+\alpha_3$};
			\node[ below, font = \tiny, black] at (3,1.5) {$\alpha_1+\alpha_2+2\alpha_3$};
			\node[ above, font = \tiny, black] at (2.4,2.4) {$\alpha_1+2\alpha_2+2\alpha_3$};
			\fill (0,0) circle (1.5pt);
			\fill (0,3) circle (1.5pt);
			\fill (3,3) circle (1.5pt);
			\fill (3,1.5) circle (1.5pt);
			\fill (2,2) circle (1.5pt);
			\fill (4,2) circle (1.5pt);
			\fill (0,6) circle (1.5pt);
			\fill (6,0) circle (1.5pt);
			\fill (2.4,2.4) circle (1.5pt);
			\node[ left, black] at (3,2) {$K$};
			\node[ left, black] at (1,1) {$K^1$};
			\node[ left, black] at (1.5,4) {$K^2$};
			\node[ left, black] at (5,1) {$K^3$};
		\end{tikzpicture}
		\caption{\centering 
			The cones $K$, $K^1$, $K^2$, $K^3$ in the case $\DB_3$.}
		\label{fig:conesB3}
	\end{figure}
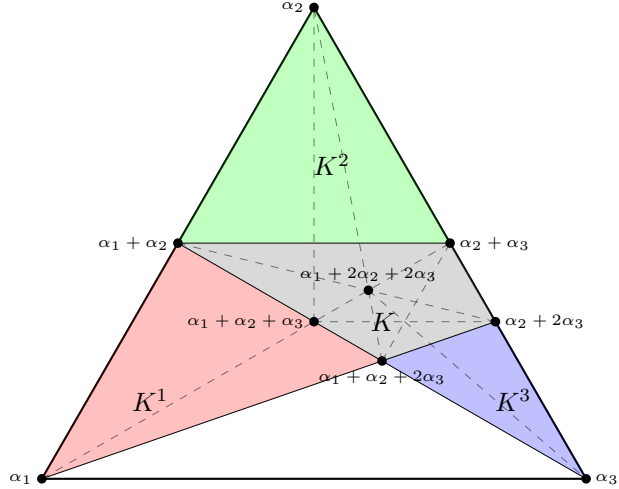

We illustrate the cones $K,K^1,K^2,K^3$ in the $B_3$ case in Figure \ref{fig:conesB3}. 

\begin{proof}
The proof follows the same line of argumentation as in Proposition \ref{prop:MUe}. 
Fix $i\in\{1,\dots,n\}$, and $u\in\Int(K^i)$; note that the region $K^i$ is a cone of maximal dimension, so $\Int(K^i)\neq\emptyset$.

Given $j\neq i$, we will show that the inequality $\nu_j x\leq \mu_j(u)$ cannot be eliminated from the definition of $P_u$, by proving that there exists an element $x\geq 0$ such that $\nu_\ell x \leq \mu_\ell(u)$ for every $\ell\neq j$ and $\nu_j x > \mu_j(u)$. 

Note that we have $\mu_j^{r_j}(u)> 0$, that is 
\begin{equation}\label{eq:inequ2}
\mu_j(u)<-\sum_{k\neq j}\langle \alpha_k,\alpha_j\rangle \mu_k(u).
\end{equation}
Consider now the set $S_j:=\{k\mid \langle \alpha_k,\alpha_j \rangle<0\}$ and observe that, for every $k\in S_j$,  $r_j(\alpha_k)=\alpha_k-\langle\alpha_k,\alpha_j\rangle\alpha_j$ is a positive non-simple root, corresponding to a column of the matrix $M$. 
Let $x\in \R^r$ be defined by having coordinate $\mu_k(u)$ in the position corresponding to $r_j(\alpha_k)$, for every $k\in S_j$, and zero in the remaining positions. Then, by construction we have
\[Mx=\sum_{k\in S_j}\mu_k(u)r_j(\alpha_{k}),\] 
so 
\[
\nu_\ell x=\sum_{k\in S_j}\mu_k(u)\mu_\ell(r_j(\alpha_{k}))=\begin{cases}
\mu_\ell(u), &\mbox{if }\ell\in S_j\\
0\leq\mu_\ell(u), &\mbox{if }\ell\notin S_j, \,\,\ell\neq j\\
-\sum_{k\in S_j}\mu_k(u)\langle \alpha_k,\alpha_j \rangle >\mu_j(u) &\mbox{if }\ell=j
\end{cases}
\]
where the last inequality follows from (\ref{eq:inequ2}). This concludes the proof. 
\end{proof}

Choose a general $u\in\Int(K^i)\cap \Mo(H)_\Q$, set $\MU_e^i:=\GU_{e,u}$, and let $\Sigma_u$ be the fan in $\No(T/H)_{\R}$ defining it. This variety is, in principle, defined up to SQM. The next statement shows that the minimal cone $\tau_i\in \Sigma_u$ containing the ray generated by $-\nu_i$ in its relative interior (which defines the divisor $D_i^{r_i}$ in $\MU_e$) is independent of the choice of $u$.  

\begin{lemma}\label{lem:center}
	For every $i\in\{1,\dots,n\}$ there exists a unique cone $\tau_i\subset\No(T/H)_\R$ that belongs to the fan $\Sigma_u$ for every $u\in\Int(K^i)\cap \Mo(H)_\Q$, and that contains $-\nu_i$ in its relative interior.
\end{lemma}

\begin{proof} Fix $i\in\{1,\dots,n\}$ and consider again the sets $S_i=\{j\mid \langle\alpha_j,\alpha_i\rangle< 0\}$. By definition, we have a relation
\[-\nu_i=-\sum_{\alpha\in \Phi^+\setminus\Pi}\mu_i(\alpha)e_\alpha.\] 
By (\ref{eq:muiri}) we also have that, for every $\alpha\in \Phi^+\setminus\Pi$, it holds
 \[-\mu_i(\alpha)=\mu_i^{r_i}(\alpha)+\sum_{j\in S_i}\langle\alpha_j,\alpha_i\rangle \mu_j(\alpha),\quad\mbox{and}\quad\mu_i^{r_i}(\alpha)\geq 0.\]
Then we may write
\[\begin{array}{rl}
-\nu_i&=\displaystyle\sum_{\alpha\in\Phi^+\setminus\Pi}\left(\mu_i^{r_i}(\alpha)+\sum_{j\in S_i}\langle\alpha_j,\alpha_i\rangle\mu_j(\alpha)\right)e_\alpha\\[10pt]
&=\displaystyle\sum_{\alpha\in\Phi^+\setminus\Pi}\mu_i^{r_i}(\alpha) e_\alpha+\sum_{j\in S_i}\sum_{\alpha\in\Phi^+\setminus\Pi}\langle\alpha_j,\alpha_i\rangle\mu_j(\alpha)e_\alpha\\[10pt]
&=\displaystyle\sum_{\alpha\in\Phi^+\setminus\Pi}\mu_i^{r_i}(\alpha) e_\alpha+\sum_{j\in S_i}(-\langle\alpha_j,\alpha_i\rangle)(-\nu_j).
\end{array}
\]

This tells us in particular that $-\nu_i$ lies in the relative interior of the cone $\tau_i$ generated by the set
\[
\{-\nu_j\mid j\in S_i\}\cup\{e_\alpha\mid \mu_i^{r_i}(\alpha)>0\}=\{-\nu_j\mid j\in S_i\}\cup\{e_\alpha\mid \alpha\notin \ol{K^i}\}.
\]
The proof will be finished by proving that $\tau_i$ is a cone of the fan $\Sigma_u$ defining the toric variety $\GU_{e,u}$, for every choice of $u\in\Int(K^i)\cap \Mo(H)_\Q$. Let us denote by $\tau_i'$ the face of $(\sigma^+)^\vee$ generated by the $1$-parameter subgroups $\fg_{-\alpha}\setminus\{0\}$ of $T$ with \[\alpha\in\{\alpha_j\mid j\in S_i\}\cup\{\alpha\in\Phi^+\setminus\Pi\mid \alpha\notin \ol{K^i}\}.\]  
Note that $\tau_i$ is the image of $\tau'_i$ into $\No(T/H)_\R$, hence it is enough to prove that $\tau'_i$ belongs to the set of $u$-semistable points, for every choice of $u\in\Int(K^i)\cap \Mo(H)_\Q$. But, given such a $u$, a direct reformulation of \cite[Corollary~14.2.22]{CLS} in the language of fans tells us that the corresponding set of semistable points in $U_e$ is the open subset defined by
\[\Sigma((\sigma^+)^\vee)^{\ss}_u=\{\tau\preceq (\sigma^+)^\vee\mid u \in \pi(\sigma^+\cap \tau^\perp)\}.\] 
Then we conclude by noting that $\tau'_i\in \Sigma((\sigma^+)^\vee)^{\ss}_u$, since $u$ belongs to the relative interior of the cone generated by $\{\alpha_j\mid j\notin S_i\}\cup\{\alpha\in\Phi^+\setminus\Pi\mid \alpha \in \ol{K^i}\}$.\end{proof}

As a consequence, we get the following immediate description of the rational contraction associated to the wall crossing defined by $\mu_i^{r_i}(m)=0$.

\begin{corollary}\label{cor:fundblup}
	With the notation of Lemma \ref{lem:center}, for any choice of general elements $v\in \Int(K)\cap \Mo(H)_\Q $ and $u\in\Int(K^i)\cap \Mo(H)_\Q$, the divisorial part of the rational contraction $\MU_e:=\GU_{e,v}\dashrightarrow \MU_e^i:=\GU_{e,u}$ is a (normalized) toric blowup of the $(T/H)$-orbit closure $V_i\subset\MU_e^i$ associated to the cone $\tau_i$.
	\qed
\end{corollary}

\subsection{Recipe for the MG quotient}

We now describe how to obtain, from an MG quotient of an affine patch, an MG quotient of the complete flag $X$. More precisely, in this section we will show the following: 

\begin{theorem}\label{thm:mellow}
There exists $m\in\Q_{>0}(\omega_1,\dots,\omega_n)$ such that \[\GU_{e,m}\] is an MG quotient of $U_e$ and such that an MG quotient of $X$ can be obtained from $\GU_{e,m}$ by means of $n$ 
divisorial extractions, without prescribed order. 
\end{theorem}

\subsubsection{Mellow weights}\label{ssec:appropriate}

Let us consider the following linear map
\[\kappa:\Mo(H)_\R\to \Mo(H)_\R, \quad \kappa(\omega):=\sum_{i=1}^n \langle \omega,\alpha_i\rangle \alpha_i.\]
By Remark \ref{rem:fundcochar}, its matrix with respect to the basis of fundamental weights is the transpose $C^t$ of the Cartan matrix of $G$. Thus $\kappa$ is a linear isomorphism and, in particular, an open map. 

\begin{definition}\label{def:mellow}
We say that $\omega\in\Mo(H)_\Q$ is a {\em mellow weight} if $\kappa(\omega)$ belongs to $\R_{\geq0}(\Phi^+\setminus\Pi)$ and lies outside every wall of the GIT chamber decomposition of the cone $\pi_\R(\sigma^+)$. In other words, $\omega$ is mellow if $\kappa(\omega)$ defines an MG quotient of $U_e$. 
\end{definition}

\begin{remark}\label{rem:mellow1}
By the proof of Proposition \ref{prop:MUe}, the weight $\omega$ is mellow if $\kappa(\omega)$ is a general element in the open set defined by
\[
\{x\in\Mo(H)_\R\mid \mu_i(x)>0, \mu_i^{r_i}(x)>0,\mbox{ for every }i\}.
\]
Choosing $\omega\in\Q_{>0}(\omega_1,\dots,\omega_n)$, we automatically have $\mu_i(\kappa(\omega))=\langle\omega,\alpha_i\rangle>0$, so the above conditions amount just to $
\mu_i^{r_i}(\kappa(\omega))>0$ or, equivalently
\begin{equation}\label{eq:mellow}
\langle \omega,\alpha_i\rangle>\langle \kappa(\omega),\alpha_i\rangle, \mbox{ for every }i.
\end{equation}
\end{remark}

\begin{figure}[ht]
	\centering    
		\begin{tikzpicture}
			\fill[gray, opacity = 0.1] (-1.24,-3) -- (3,-3) -- (3,1.24)-- (-5.48,1.24) -- cycle;
			\draw[thick] (-1.24,-3) -- (-5.48,1.24);
			\draw[thick] (-1.24,-3) -- (3,-3);
			\node[left, font = \footnotesize] at (-3.5,-0.8) {$\Eff(\MU_e)$};
			
			\fill[gray, opacity = 0.3] (-1.24,-3) -- (3,1.24)-- (-1.24,1.24) -- cycle;
			\node[font = \footnotesize] at (0,0.6) {$\Mov(\MU_e)$};
			\draw[thick, red] (-1.24,-3) -- (-1.24,1.24);
			\draw[thin, dashed, opacity = 0.3, red] (-3,-1.24) -- (-3,1.24);
			\draw[thin, dashed, opacity = 0.3, red] (1.24,-3) -- (1.24,1.24);
			\draw[thick, dashed, red] (-3,-1.24) -- (3,-1.24);
			\draw[thick, blue] (-1.24,-3) -- (3,1.24);
			\draw[thin, dashed, opacity = 0.3, blue] (-3,-1.24) -- (-0.51,1.24);
			\draw[thin, dashed, opacity = 0.3, blue] (3,-1.24) -- (0.51,1.24);
			\draw[thin, dashed, opacity = 0.3, blue] (1.24,-3) -- (3,-1.24);
			\draw[thick, dashed, blue] (1.24,-3) -- (-3,1.24);
			
			\fill (-0.51,-1.24) circle[radius = 1pt];
			\node[below left, font = \tiny] at (-0.3,-1.24) {$\kappa(\omega)$};
			
			\fill (0,0) circle[radius = 1pt];
			\node[below, font = \tiny] at (0,0) {$\omega$};

			\draw[->,  thick, gray](-1.24,-3) -- (-0.24,-3);
			\draw[->,  thick, gray](-1.24,-3) -- (-2.24,-2);
			\draw[->,  thick, gray](-1.24,-3) -- (-1.24,-2);
			\draw[->,  thick, gray](-1.24,-3) -- (-0.24,-2);
			\fill (-1.24,-3) circle[radius = 1.5pt];   
			\node[ left, font = \tiny] at (-1.25,-3) {$eB$};
						
		\end{tikzpicture}
\caption{\centering The cones of divisors of $\MU_e$ and a mellow weight $\omega$.}
\label{fig:cones}
\end{figure}
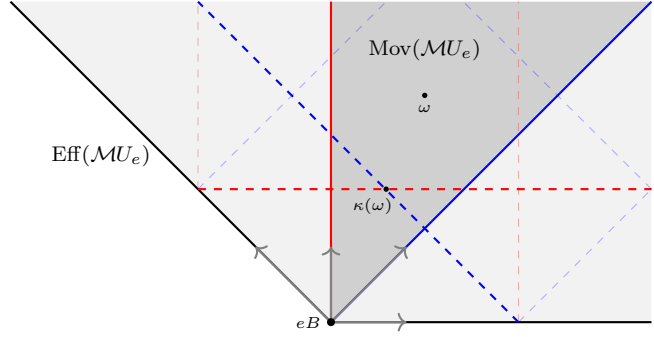

We may intuitively summarize this by means of Figure \ref{fig:cones}, where we illustrate, in the type $\DC_2$ case, the effective and movable cones of $\MU_e$. Note that the weight $\omega$ is mellow since $\kappa(\omega)$ lies in the movable cone's interior.  

\begin{proposition}\label{prop:existmellow}
There exists a mellow weight $\omega\in\Q_{>0}(\omega_1,\dots,\omega_n)$.
\end{proposition}

In order to prove this statement we will need some preliminary results. Let us introduce the {\em Cartan adjacency matrix} of $G$, defined as the non-negative matrix
\[\Theta:=2I-C^t.\]
This matrix defines a strongly connected directed graph whose underlying undirected graph is the Dynkin diagram of $G$. The Perron-Frobenius theorem asserts that the spectral radius $\rho(\Theta)$ is a positive real eigenvalue of $\Theta$ with geometric multiplicity equal to one. Let us denote by $\theta$ an associated eigenvector of $\rho(\Theta)$, which is usually called a {\em Perron vector}. 

\begin{lemma}\label{lem:spectralradius}
With the above notation, if the Dynkin diagram of $G$ is different from $\DA_1$, $\DA_2$ (in particular, in the situation of Setup \ref{set:flags}), then $\rho(\Theta)>1$. 
\end{lemma}

\begin{proof}
Given a subset $I\subset\{1,\dots,n\}$, let us denote by $\Theta(I)$ the $n\times n$ matrix with coefficients
\[
\Theta(I)_{ij}=\begin{cases}\Theta_{ij},&\mbox{if }i,j\in I\\
0,&\mbox{otherwise.}\end{cases}
\]
Let us also denote by $\Theta(I)'$ the $|I|\times |I|$ matrix constructed by eliminating in $\Theta(I)$ all the zero columns and rows.
In particular, we have $\Theta_{ij}\geq \Theta(I)_{ij}\geq 0$ for every $i,j$ and, subsequently, $\Theta^k_{ij}\geq \Theta(I)^k_{ij}\geq 0$ for $k>0$. It then follows that 
\[\|\Theta^k\|^{1/k}\geq \|\Theta(I)^k\|^{1/k},\] 
where $\|\cdot\|$ denotes the $\infty$-norm. By Gelfand's formula, we find that \[\rho(\Theta)\geq \rho(\Theta(I))=\rho(\Theta(I)'), \quad\textrm{for every }I.\] 

We conclude by noting that, in the situation of the statement, we may assume that there exists $I$ such that $\Theta(I)'$ is the 
Cartan adjacency matrix associated to a group of type $\DA_3$, $\DB_2$ or $\DG_2$, whose spectral radii are $\sqrt{2},\sqrt{2},\sqrt{3}$, respectively.
\end{proof}

\begin{proof}[Proof of Proposition \ref{prop:existmellow}]
Given a positive Perron vector $\theta=(a_1,\dots,a_n)\in\R^n$, we can define the element
\[\omega_\theta:=\sum_{i=1}^n a_i\omega_i\in\Mo(H)_\R.\]
For every $i\in\{1,\dots,n\}$, we have
\[\begin{array}{rl}
\langle \kappa(\omega_\theta),\alpha_i\rangle&=\langle\sum_ja_j\alpha_j,\alpha_i\rangle=\sum_ja_j\langle\alpha_j,\alpha_i\rangle=\sum_ja_jC_{ji}=\\[4pt]
&=(2-\rho(\Theta))a_i=(2-\rho(\Theta))\langle \omega_\theta,\alpha_i\rangle<\langle \omega_\theta,\alpha_i\rangle,
\end{array}
\]
where the last inequality follows from Lemma \ref{lem:spectralradius}. On the other hand, since all the eigenvalues of the Cartan matrix are positive, we also have $2-\rho(\Theta)>0$, which shows that $\kappa(\omega_\theta)$ is strictly dominant. Both properties hold in a sufficiently small open neighborhood of $\omega_\theta$. Choosing a rational weight $\omega$ in this neighborhood such that $\kappa(\omega)$ avoids the walls of the GIT chamber decomposition, we conclude that $\omega$ is mellow by Remark~\ref{rem:mellow1} and that $\kappa(\omega)$ is strictly dominant.
\end{proof}

\begin{proof}[Proof of Theorem \ref{thm:mellow}]
We start by fixing a mellow weight $\omega$ and consider the MG quotient 
\begin{equation*}\label{eq:mellowMG}
\MU_e:=\GU_{e,\kappa(\omega)}.
\end{equation*}

Let us observe that, by definition of $\kappa$, the point $\kappa(\omega)$ is the intersection of the $n$ hyperplanes 
\[
\{v\in\Mo(H)_\R\mid \mu_i(v)=\mu_i(\omega)-\mu_i^{r_i}(\omega)\}, \qquad i=1,\ldots,n.
\]
These hyperplanes correspond to the GIT walls associated to the fundamental divisors $D_i^e$. Since the cocharacters $\mu_1,\ldots,\mu_n$ are linearly independent, they intersect transversally at $\kappa(\omega)$. Moreover, since $\omega$ is mellow, no other GIT walls meet a sufficiently small neighborhood of this point. Consequently, the local GIT chamber decomposition around $\kappa(\omega)$ is equivalent to the decomposition of $\R^n$ by its coordinate hyperplanes.

We can therefore pass from a chamber corresponding to an MG quotient of $U_e$ to a chamber contained in $\Delta^{\mov}(\omega)$ by crossing the $n$ hyperplanes above in any prescribed order, without crossing any additional walls. Since each of these wall crossings corresponds to the extraction of the fundamental divisor $D_i^e$, we conclude that an MG quotient of $X$ can be obtained from $\MU_e$ by performing the corresponding $n$ divisorial extractions in any order. 
\end{proof}

\subsubsection{Divisorial extractions and $\MX$}\label{ssec:divextr}

Thanks to the above arguments, we know that there exists a sequence of $n$ divisorial extractions leading from $\MU_e$ to $\MX$, which can be performed in any order. As noticed in Remark \ref{rem:MXuponMU}, the exceptional divisors will be $D_1^e,\ldots,D_n^e$; we conclude by  describing such divisorial extractions. 

For each $i=1,\ldots,n$, set
\[\theta_i:=\operatorname{ord}_{D_i^{r_i}}\]
be the toric divisorial valuation associated with the fundamental contraction of Corollary \ref{cor:fundblup}. The Weyl group acts on the GIT chamber decomposition of the weight polytope of $X$: the action of $w\in W$ sends the semistable locus for a linearization $u$ onto the semistable locus for the transformed linearization $w(u)$, and hence identifies the corresponding wall crossings. In particular, the reflection $r_i$ sends the wall crossing associated with $D_i^{r_i}$ to the one associated with $D_i^e$, whose corresponding divisorial valuation on $\C(\MU_e)$ will be:
\[
\nu_i:=r_i^*\theta_i. 
\]

The induced birational involution $r_i:\MU_e\dashrightarrow\MU_e$ is obtained from the Weyl action $U_e\to U_{r_i}$ and the rational identification $U_{r_i}\dashrightarrow U_e$ defined by their intersection, after passing to the chosen GIT quotients. Since $r_i$ maps $\overline{X^-(Y_i^{r_i})}$ onto $\overline{X^+(Y_i^e)}$, the image of the generic point of $D_i^{r_i}$ parametrizes an $H$-orbit in the latter locus. Consequently the center of $\nu_i$ on $\MU_e$ is
\[Z_i:=\overline{r_i(D_i^{r_i})}.\]
Equivalently, $Z_i$ is the closure of the corresponding locus of $H$-orbits in $\overline{X^+(Y_i^e)}$.

Moreover, the transported wall crossing is divisorial. Denote by $f_i:M_i\to\MU_e$ its associated extraction and by $E_i$ its exceptional divisor, so that $\operatorname{ord}_{E_i}=\nu_i$. Since $-E_i$ is relatively ample (after taking a suitable positive multiple), the valuation ideal algebra
\[\mathcal R(\nu_i):=\bigoplus_{m\geq0}\mathfrak a_m(\nu_i),\qquad\mathfrak a_m(\nu_i):=\{f\in\mathcal O_{\MU_e}\mid\nu_i(f)\geq m\},\]
is finitely generated, and the corresponding divisorial extraction is the normalization of 
\[\operatorname{Proj}_{\MU_e}\mathcal R(\nu_i). \]
Thus the $n$ divisorial extractions required to construct $\MX$ are the blowups of the valuations associated with $r_1^*\theta_1,\ldots,r_n^*\theta_n$ that, as established in the proof of Theorem~\ref{thm:mellow}, can be performed in any prescribed order.

\bibliographystyle{plain}
\bibliography{bibliomin}

\end{document}